\documentclass[12pt]{amsart}
\usepackage{amssymb}
\usepackage{amsmath}
\usepackage{mathabx}
\usepackage{amsthm}
\usepackage{amstext}
\usepackage{amsopn}
\usepackage{mathrsfs}
\usepackage{latexsym}
\DeclareMathAlphabet{\mathpzc}{OT1}{pzc}{m}{it}
\allowdisplaybreaks
\newtheorem{Definition}{Definition}[section]
\newtheorem{Proposition}{Proposition}[section]
\newtheorem{Lemma}{Lemma}[section]
\newtheorem{Theorem}{Theorem}[section]
\newtheorem{Corollary}{Corollary}[section]
\newtheorem{Remark}{Remark}[section]
\newtheorem{Example}{Example}[section]

\begin{document}
%%%%%%%%%%%%%%%%%%%%%%%%%%%%%%%%%%%%%%%%%%%%%%%%%%%%%%%%%%%%
\bibliographystyle{plain}
\footnotetext{
This work is supported by Narodowe Centrum Nauki grant No. 2014/15/B/ST1/00166}
%%%%%%%%%%%%%%%%%%%%%%%%%%%%%%%%%%%%%%%%%%%%%%%%%%%%%%%%%%%%%%%%%
\title[Random matrices and continuous circular systems]
{Random matrices, continuous circular systems\\ and the triangular operator}
\author[R. Lenczewski]{Romuald Lenczewski}
\address{Romuald Lenczewski, \newline
Wydzia\l{} Matematyki, Politechnika Wroc\l{}awska, \newline
Wybrze\.{z}e Wyspia\'{n}skiego 27, 50-370 Wroc{\l}aw, Poland  \vspace{10pt}}
\email{Romuald.Lenczewski@pwr.edu.pl}
%%%%%%%%%%%%%%%%%%%%%%%%%%%%%%%%%%%%%%%%%%%%%%%%%%%%%%%%%%%%%%%%
\begin{abstract}
Using suitably defined continuous analogs of the matricial circular systems 
and the direct integral of Hilbert spaces ${\mathcal H}=\int_{\Gamma}^{\oplus}{\mathcal H}(\gamma)d\gamma$, 
we study the operators living in ${\mathcal H}$ which give the asymptotic joint *-distributions of 
complex independent Gaussian random matrices with not necessarily equal variances of the entries. These operators are decomposed in terms of {\it continuous circular systems}
$\{\zeta(x,y;u):x,y\in [0,1], u\in \mathpzc{U}\}$ acting between the fibers of ${\mathcal H}$,
the continuous analogs of matricial circular systems 
obtained when the Gaussian entries are block-identically distributed.
In the case of square matrices with i.i.d. entries, we obtain the circular operators of Voiculescu, 
whereas in the case of upper-triangular matrices with i.i.d. entries, we obtain the triangular operators of Dykema and Haagerup.
We apply this approach to give a bijective proof of the formula for the moments of $T^*T$, where $T$ is a triangular operator, using the enumeration formula of Chauve, Dulucq and Rechnizter for alternating ordered rooted trees.
\end{abstract}

\maketitle

\section{Introduction}
Independent Gaussian random matrices with suitably normalized complex i.i.d. entries are asymptotically free with respect to the normalized trace composed with classical expectation. The limit mixed (*-) moments can be expressed in terms of mixed (*-) moments of free circular operators. This fundamental result was shown by Voiculescu [16], who also found a relation to free group factors [17].

Asymptotic freeness of Voiculescu was later generalized in many directions, 
in particular, to the non-Gaussian random matrices by Dykema [3] and to Gaussian band matrices by Shlyakhtenko [13]. In the latter case, where the Gaussian variables are not assumed to be identically distributed, scalar valued freeness [18] is not sufficient to
describe the asymptotics of matrices and one has to use freeness with amalgamation over some subalgebra, a generalization of freeness, in which a state is replaced by an operator-valued conditional expectation  
with values in this subalgebra. This approach was later used by Benaych-Georges [1] in his study 
of rectangular block random matrices, who introduced a rectangular analog of Voiculescu's R-transform [18].
Analytic methods like operator-valued transforms were also applied to the study of blocks of random matrices (most of these results are mentioned in the recently published monograph of Mingo and Speicher [11]). 

The main point of our approach is that our realizations of the limit joint *-distributions of Gaussian random matrices are built from
operators living in Hilbert spaces rather than in Hilbert modules. 
Our combinatorics also has new features since it is based on coloring
noncrossing pair partitions rather than on using nested evaluations which are needed when applying the conditional expectation 
and operator-valued free probability. Although one can translate our approach to that of the operator-valued free probability, the new language provides a new realization of limit distributions as well as connections to some nice results in combinatorics. For instance, 
the combinatorics of colored partitions allows us to find a connection between 
certain results on moments of the triangular operator defined by Dykema and Haagerup [4] with
known enumeration results for alternating ordered rooted trees. The triangular operator 
got some attention in the context of the famous invariant subspace problem in the papers of 
Dykema and Haagerup on the so-called DT-operators and their decompositions [4,5], 
including the circular operator of Voiculescu and the triangular operator $T$. The moments of $T^*T$
were also found by Dykema and Haagerup in [4] and the more general case of the moments of ${T^*}^kT^k$ was 
treated by \'Sniady [14], who used the 
so-called generalized circular elements and the combinatorics of nested evaluations.

Our original approach to the study of block matrices was related directly to 
the connection between free probability and random matrices [7,8,9,10]. 
Since independent (Gaussian and also more general)
random matrices $\{Y(u,n): u\in \mathpzc{U}\}$, are asymptotically free under the expectation of 
the trace as $n\rightarrow \infty$, their blocks 
$$
\{S_{p,q}(u): 1\leq p,q\leq r, u\in \mathpzc{U}\}
$$
should be asymptotically 'matricially free' under `partial traces'. This idea was actually 
our main motivation to introduce a concept 
of a `matricial' form of independence, called {\it matricial freeness} [7], where
a family of subalgebras of a given algebra is replaced by an array of its subalgebras and a state is replaced by a family of states corresponding to the array's diagonal. 
Although our construction was entirely based on scalar-valued states, 
it turned out later that this concept was related to freeness with amalgamation and operator-valued 
conditional expectation, at least in the context of {\it matricial circular systems} studied in [10], where we refer the reader for details on this relation. 

As far as matricial circular systems are concerned, they describe the asymptotics of blocks of independent Gaussian 
random matrices under partial traces.
If we are given an ensemble of independent non-Hermitian $n\times n$ 
Gaussian random matrices $\{Y(u,n):u\in \mathpzc{U}\}$ with suitably normalized 
independent block-identically distributed ({\it i.b.i.d.}) complex entries for each natural $n$, 
then the mixed *-moments of their (in general, rectangular) blocks $S_{p,q}(u,n)$ 
converge under normalized partial traces to the mixed 
*-moments of certain bounded operators, which we write informally
$$
\lim_{n\rightarrow \infty}S_{p,q}(u,n)=\zeta_{p,q}(u),
$$
where $u\in \mathpzc{U}$ and $1\leq p, q \leq r$ and
the operators $\zeta_{p,q}(u)$ are called 
{\it matricial circular operators} [10].
The arrays of such operators play the role of matricial analogs of 
circular operators. If we want to find the  
counterparts of these operators in the operator-valued free probability, 
we need to take operators of the form $F_p c(u) F_q$, where $\{c(u):u\in \mathpzc{U}\}$ is a family
of circular elements living in the Fock space over Hilbert $A$-module, where $A$ is the algebra of $r\times r$ diagonal matrices with canonical generators $F_1, \ldots , F_r$.

In particular, if the matrices $Y(u,n)$ are square and have i.i.d. 
(normalized in a standard way) Gaussian entries and we divide them into $r^2$ rectangular 
blocks $S_{p,q}(u)$, whose sizes tend to
$d_{p}\times d_q$ as $n\rightarrow \infty$, respectively, where $d_1+\ldots +d_r=1$, then the limit joint *-distribution 
of the family $\{Y(u,n):u\in \mathpzc{U}\}$ under the expectation of the normalized 
trace is that of the family of standard circular operators $\{\zeta(u):u\in \mathpzc{U}\}$ and these are 
decomposed in terms of matricial circular operators as
$$
\zeta(u)=\sum_{p,q=1}^{r}\zeta_{p,q}(u)
$$
for any $u\in\mathpzc{U}$.
They are expressed in terms of the canonical matricial creation and annihilation operators as 
$$
\zeta_{p,q}(u)=\wp_{p,q}(u')+\wp_{q,p}^*(u''),
$$
which corresponds to the the well-known realization of the 
circular operator of Voiculescu $c=\ell_1+\ell_2^*$, where 
$\ell_1, \ell_2$ are free. 
Here, we use the indices $u',u''$ to encode the fact that 
we have two free copies of creation-annihilation pairs 
related to the array labeled by $u$. 
If matrices have entries 
which are only i.b.i.d. 
and not i.i.d., then a similar formula holds except that operators 
$\zeta(u)$ are no longer circular.
In particular, this includes the case when the matrices $Y(u,n)$ are block strictly upper triangular 
(it suffices to set the covariances of all variables lying on and under the diagonal to be equal to zero). 
In a similar manner, block lower triangular matrices can be considered. 

One has to stress that objects like matricial circular operators are not new. 
Voiculescu introduced and studied such objects in [17], although perhaps not in
the context of the asymptotics of rectangular blocks.
Namely, if ${\mathcal A}$ is the $C^{*}$-algebra generated by the system of 
free (generalized) circular operators
$$
\{c(p,q,u): 1\leq p, q \leq r, u\in \mathpzc{U}\},
$$
then
$$
\zeta_{p,q}(u)=c(p,q,u)\otimes e(p,q)\in M_{r}({\mathcal A})
$$
for any $p,q,u$. These operators are the generators of free group factors used by 
Voiculescu in his proof of free group factors isomorphisms [17, Theorem 3.3]. 
We defined them differently, as operators on the matricially free Fock space of tracial type, but 
we showed in [10] that our matricial circular operators were in fact of the same tensor product form.  

In this paper, we generalize the decomposition of the operators 
$\zeta(u)$, which give the asymptotic *-joint distributions of matrices $Y(u,n)$, to the continuous case. Namely, we investigate the limits of the *-mixed moments obtained in the discrete case as $r\rightarrow \infty$
and show that they agree with the *-mixed moments of operators that can be written in the form
$$
\zeta(g,u)=\int_{\Gamma_1}^{\oplus}g(x,y)\zeta(x,y;u)dxdy,
$$
where $\Gamma_1=[0,1]\times [0,1]$, $g\in L^{\infty}(\Gamma_1)$, 
and the operators $\zeta(x,y;u)$ are continuous generalizations of the matricial circular operators
$\zeta_{p,q}(u)$: roughly speaking, $p/r\rightarrow x$, $q/r\rightarrow y$ as $r\rightarrow \infty$.
The operators $\zeta(x,y;u)$ are continuous analogs of the matricial circular operators $\zeta_{p,q}(u)$, namely
$$
\zeta(x,y;u)=\wp(x,y;u')+\wp^*(y,x;u'')
$$
for any $x,y\in [0,1]$ and $u\in \mathpzc{U}$, where 
$$
\{\wp(x,y;u):x,y\in [0,1], u\in \mathpzc{U}'\cup\mathpzc{U}''\}
$$
is a suitably defined continuous family of creation operators, where $\mathpzc{U}'$ and $\mathpzc{U}''$ denote
replicas of $\mathpzc{U}$ and the star stands for the adjoint operation. 

In order to define all these operators, we introduce a continuous analog of the 
{\it matricially free Fock space of tracial type}, namely the direct integral of Hilbert spaces 
$$
{\mathcal H}=\int_{\Gamma}^{\oplus} {\mathcal H}(\gamma)d\gamma ,
$$
where $\Gamma=\bigcup_{n=0}^{\infty}\Gamma_{n}$ and 
$\Gamma_{n}=I^{n+1}$, with $I=[0,1]$ ($d\gamma$ stands for 
the direct sum of Lebesgue measures on various $\Gamma_n$), equipped 
with the state
$$
\varphi(F)=\int_{I} \langle F(x)\Omega(x), \Omega(x)\rangle dx,
$$
where $dx$ is the Lebesgue measure on $I$ and 
$F=\int_{I}^{\oplus}F(x)dx$ is a decomposable operator field.
In order to decompose the operators  $\zeta(g,u)$ (circular, triangular, or other) as integrals over $\Gamma_1$
(or its subset, if $g$ vanishes outside this subset) one uses various decompositions of ${\mathcal H}$ into direct integrals.

In particular, if $g=\chi_{\Delta}$, where $\Delta=\{(x,y)\in \Gamma_1: x<y\}$, 
then we obtain the integrals 
$$
T(u)=\int_{\Delta}^{\oplus}\zeta(x,y;u)dxdy ,
$$ 
which play the role of continuous decompositions of the family of free triangular operators 
of Dykema and Haagerup [4]. We apply this approach to study the mixed *-moments of the triangular operators. 
In particular, we provide a new bijective proof of the formula for the moments 
$$
\varphi((T^*T)^{n})=\frac{n^n}{(n+1)!},
$$
shown by Dykema and Haagerup [4] by a different method, where $T$ is a triangular operator, using the nice enumeration result of Chauve, Dulucq and Rechnitzer [2] for alternating ordered rooted tress. In this context, let us remark that there is a more general bijective proof of \'Sniady [15], based on a certain 
algorithm of counting total orders on directed trees. 

The paper is organized as follows. In Section 2 we recall the notions related to matricial circular systems.
In Section 3, we introduce the direct integral of Hilbert spaces, on which we define a family
of creation and annihilation operators. In Section 4, we introduce continuous circular systems as 
isometries between suitably defined fiber Hilbert spaces. Mixed *-moments of 
matricial circular operators and their convergence to the mixed *-moments of the operators  
decomposed in terms of the continuous circular systems are discussed in Section 5. 
In Section 6, we apply our approach to the triangular operators and we provide a bijective proof 
of the formula for the moments of $T^*T$.

We adopt the convention that the stars which indicate adjoints are written closely to the main symbol, for instance
$\wp^*(g,u)$ and $\wp^*(x,y;u)$ are the adjoints of $\wp(g,u)$ and $\wp(x,y;u)$, respectively.

\section{Matricial circular systems}

The matricial circular operators studied in [10] live in a kind of matricially free Fock space.
The original definition given in [7] was slightly generalized in [9,10], where we used the name 
of the `matricially free Fock space of tracial type' which will be used below.

Let $[r]:=\{1,2, \ldots, r\}$ and let $\mathpzc{U}$ be a countable set.
To each $(p,q)\in \mathcal{J}\subset [r]\times [r]$ and $u\in \mathpzc{U}$ we associate a Hilbert space 
${\mathcal H}_{p,q}(u)$. Using this family of Hilbert spaces, we can construct
the matricially free Fock space of tracial type.

\begin{Definition}
{\rm By the {\it matricially free Fock space of tracial type} we understand the direct sum of Hilbert spaces
\begin{equation*}
{\mathcal M}= \bigoplus_{q=1}^{r} {\mathcal M}_{q},
\end{equation*}
where each summand is of the form
\begin{equation*}
{\mathcal M}_{q}={\mathbb C}\Omega_{q}\oplus \bigoplus_{m=1}^{\infty}
\bigoplus_{\stackrel{p_1,\ldots , p_m}
{\scriptscriptstyle u_1, \ldots , u_n}
}
{\mathcal H}_{p_1,p_2}(u_{1})\otimes {\mathcal H}_{p_2,p_3}(u_2)
\otimes \ldots \otimes 
{\mathcal H}_{p_m,q}(u_{m}),
\end{equation*}
where $\Omega_q$ is a unit vector for any $q\in [r]$, endowed with the canonical inner product. 
We denote by $\Psi_q$ the state associated with $\Omega_q$.}
\end{Definition}

Let us recall a number of basic facts and notions from [8,9,10].
\begin{enumerate}
\item 
In the special case when each ${\mathcal H}_{p,q}(u)={\mathbb C}e_{p,q}(u)$ for any $p,q,u$, where $e_{p,q}(u)$ is a unit vector, the canonical orthonormal basis of ${\mathcal M}$ consists of tensors of the form
$$
e_{p_1,p_2}(u_1)\otimes e_{p_2,p_3}(u_2)\otimes  \ldots \otimes e_{p_m,q}(u_{m}),
$$
where $p_1, \ldots, p_m,q\in [r]$, $u_{1},\ldots , u_m\in \mathpzc{U}$ and $m\in {\mathbb N}$, and of vacuum vectors $\Omega_1, \ldots , \Omega_r$. In the general case, a canonical basis of similar form can be given, except that the basis of each ${\mathcal H}_{p,q}(u)$ may contain more than one vector.  
\item
The {\it matricially free creation operators} associated with matrices 
$B(u)=(b_{p,q}(u))$ of non-negative real numbers (covariance matrices) are bounded
operators whose non-trivial action onto the basis vectors is  
\begin{eqnarray*}
\wp_{p,q}(u)\Omega_q&=&\sqrt{b_{p,q}(u)}e_{p,q}(u),\\
\wp_{p,q}(u)(e_{q,t}(s))&=&\sqrt{b_{p,q}(u)}(e_{p,q}(u)\otimes e_{q,t}(s)),\\
\wp_{p,q}(u)(e_{q,t}(s)\otimes w)&=&\sqrt{b_{p,q}(u)}(e_{p,q}(u)\otimes e_{q,t}(s)\otimes w),
\end{eqnarray*}
for any $p,q,t\in [r]$ and $u,s\in \mathpzc{U}$, where $e_{q,t}(s) \otimes w$ is a basis vector.
Action onto the remaining basis vectors gives zero. The corresponding {\it matricially free annihilation operators} 
are their adjoints denoted $\wp_{p,q}^{*}(u)$. If $b_{p,q}(u)=1$, we call the associated operators {\it standard}.
\item
Matricially free creation operators can be realized as 
operator-valued matrices [10].
If we are given a $C^{*}$-probability space $({\mathcal A}, \phi)$ and a family of free creation operators
$\{\ell(p,q,u): p,q\in [r], \;u\in \mathpzc{U}\}$
with covariances $b_{p,q}(u)$, respectively, which are *-free with respect to $\phi$,
and $(e(p,q))$ is the array of matrix units in $M_{r}({\mathbb C})$, then
\begin{eqnarray*}
\wp_{p,q}(u)&=&\ell(p,q,u)\otimes e(p,q)
\end{eqnarray*}
for any $p,q\in [r]$ and $u\in \mathpzc{U}$. This equality holds
in the sense that the mixed *-moments of the operators $\wp_{p,q}(u)$ under the states $\Psi_j$ agree with the corresponding mixed 
*-moments of the above matrices under the states $\Phi_j=\phi\otimes \psi_j$, where $\psi_j$
is the state associated with the canonical basis vector $e(j)$ of ${\mathbb C}^{r}$, where $j\in [r]$.
\item
If $\ell^*(p,q,u)$ is the free annihilation operator 
corresponding to $\ell(p,q,u)$, then
$$
\ell^*(p,q,u)\ell(p',q',u')=\delta_{p,p'}\delta_{q,q'}\delta_{u,u'}b_{p,q}(u)
$$ 
for any $p,q,u,p',q',u'$.
\item
The {\it matricially free Gaussian operators} (called also {\it matricial semicircular operators}) take the form
$$
\omega_{p,q}(u)=\wp_{p,q}(u)+\wp^{*}_{p,q}(u),
$$
where $p,q\in [r]$ and $u\in \mathpzc{U}$. We give this definition for completeness since we will 
not use these operators in this paper.  
\item
The {\it matricial circular operators} are obtained from 
a family of arrays of matricially free creation operators of the form
$$
\zeta_{p,q}(u)=
\wp_{p,q}(u')+\wp_{q,p}^*(u''),
$$
where $u\in \mathpzc{U}=[t]$, $p,q\in [r]$ and $u',u''$ are different copies of $u$. 
As we showed in [9], they can be realized as matrices of the form
$$
\zeta_{p,q}(u)=c(p,q,u)\otimes e(p,q),
$$
where $p,q\in [r]$, $u\in \mathpzc{U}$,
$\{c(p,q,u):p,q\in [r], \;u\in \mathpzc{U}\}$ 
is a family of free generalized circular operators, i.e. $c(p,q,u)=\ell(p,q,u')+\ell^*(q,p,u'')$
and $(e(p,q))$ is an array of matrix units.
We use the term `generalized' since, in general, the covariances of the creation operators 
are arbitrary. If all creation operators are standard, this is a family of free circular operators.
The above equality holds in the same sense as in (3). The corresponding family of arrays of operators is 
called the {\it matricial circular system}.
\end{enumerate}

\section{Direct integrals}

We would like to construct a continuous analog of the matricially free 
Fock space of tracial type. 
For that purpose, we will use the formalism of direct integrals 
(for more on direct integrals, see, for instance, [6]).
However, two different direct integral decompositions of the considered Fock space
and of the canonical operator fields acting on this Fock space will be helpful.

We begin with a decomposition which is a straightforward generalization of 
the discrete matricially free Fock space of tracial type. 
For $I=[0,1]$, let 
$$
\Gamma=\bigoplus_{n=0}^{\infty}\Gamma_{n}
$$
be the direct sum of measure spaces, where $\Gamma_{n}=I^{n+1}$ is equipped with the Lebesgue measure denoted $d\gamma_n$, 
and let us denote by $d\gamma$ the corresponding direct sum of measures 
on the set $\Gamma$.

\begin{Definition}
{\rm By the {\it continuous matricially free Fock space} we understand
the direct integral of Hilbert spaces of the form 
$$
{\mathcal H}=\int_{\Gamma}^{\oplus}{\mathcal H}(\gamma)d\gamma,
$$
where Hilbert spaces are associated to $\gamma\in \Gamma$ as follows:
\begin{enumerate}
\item
if $\gamma=x\in \Gamma_{0}=I$, then 
$$
{\mathcal H}(\gamma)={\mathbb C}\Omega(x),
$$
where $\Omega(x)$ is a unit vector,
\item
if $\gamma=(x_1,x_2, \ldots , x_{n+1})\in \Gamma_{n}$ and $n\in \mathbb{N}$, then 
$$
{\mathcal H}(\gamma)=
{\mathcal H}(x_1,x_2)\otimes {\mathcal H}(x_2,x_3)\otimes \ldots \otimes {\mathcal H}(x_{n},x_{n+1}),
$$
where each ${\mathcal H}(x,y)$ is a separable Hilbert space, 
and each ${\mathcal H}(\gamma)$ is equipped with the canonical inner product,
\item
the canonical inner product in ${\mathcal H}$ is then given by
$$
\langle F, G \rangle =\int_{\Gamma} \langle F(\gamma), G(\gamma) \rangle \,
d\gamma ,
$$
where $F=\int_{\Gamma}^{\oplus}F(\gamma)d\gamma, 
G=\int_{\Gamma}^{\oplus}G(\gamma)d\gamma\in \mathcal{H}$ are 
measurable square integrable fields with the natural assumption that 
$F(\gamma),G(\gamma)\in {\mathcal H}(\gamma)$. 
\end{enumerate}}
\end{Definition}

\begin{Remark}
{\rm Let us collect certain basic facts about the Hilbert spaces defined above.
\begin{enumerate}
\item
The continuous family of unit vectors
$\{\Omega(x):x\in I\}$ replaces the finite set of vacuum vectors  $\{\Omega_1, \ldots , \Omega_r\}$
used in the discrete case. The corresponding direct integral
$$
\mathcal{H}_{0}:=\int_{I}^{\oplus}\mathcal{H}(x)dx\cong L^{2}(I)
$$ 
will be called the {\it vacuum space}. In this paper,
we will be mainly concerned with the function on $I$ which is constantly equal to one since 
it corresponds to the canonical trace on the algebra of random matrices. However, 
weighted traces will, in general, lead to different elements of $L^{2}(I)$.
\item
In the case when ${\mathcal H}(x,y)\cong {\mathcal G}$ for any $(x,y)\in \Gamma_1$, where
${\mathcal G}$ is a separable Hilbert space (with an orthonormal basis indexed by $\mathpzc{U}$),
we also have isomorphisms for higher order integrals
$$
\mathcal{H}_{n}:=\int_{\Gamma_{n}}^{\oplus}\mathcal{H}(\gamma)d\gamma_{n}
\cong L^{2}(\Gamma_{n},{\mathcal G}^{\otimes n}),
$$   
where $n\geq 1$ and $L^2(\Gamma_n,{\mathcal G}^{\otimes n})$ denotes the Hilbert space of square integrable 
${\mathcal G}^{\otimes n}$-valued functions over the set $\Gamma_n$ with respect to $d\gamma_n$ ($d\gamma$ restricted to $\Gamma_n$).
Thus, in this particular case, we have the isomorphism
$$
{\mathcal H}\cong 
L^{2}(I)\oplus \bigoplus_{n=1}^{\infty}L^{2}(\Gamma_{n}, \mathcal{G}^{\otimes n}).
$$
\item
Fields $F=\int_{\Gamma}^{\oplus}F(\gamma)d\gamma, 
G=\int_{\Gamma}^{\oplus}G(\gamma)d\gamma\in \mathcal{H}$ 
have direct sum decompositions
$$
F=\sum_{n=0}^{\infty}F_n\;\;{\rm and}\;\;G=\sum_{n=0}^{\infty}G_n,
$$
where $F_{n}, G_{n}\in \int_{\Gamma_{n}}{\mathcal H}(\gamma)d\gamma$ in the natural sense.
Under the above isomorphism assumptions, $F_0,G_0\in L^{2}(I)$ and 
$F_n,G_n\in L^{2}(\Gamma_n, \mathcal{G}^{\otimes n})$ for $n\geq 1$.
In most computations, it is enough to consider these to be of the form
\begin{eqnarray*}
F_{n}(\gamma)&=&f_1(x_1,x_2)\otimes \ldots \otimes f_{n}(x_{n},x_{n+1}),\\
G_{n}(\gamma)&=&g_1(x_1,x_2)\otimes \ldots \otimes g_{n}(x_{n},x_{n+1}),
\end{eqnarray*}
for $\gamma=(x_1, \ldots , x_{n+1})$ and $n\geq 1$, with
$f_i(x_i,x_{i+1}),g_i(x_i,x_{i+1})\in\mathcal{G}$ for any $i$.
\item
The canonical inner product in  ${\mathcal H}$ decomposes as 
$$
\langle F, G \rangle =\sum_{n=0}^{\infty}\int_{\Gamma_{n}} \langle F_n(\gamma), G_n(\gamma) \rangle \,d\gamma_n
$$
for any $F,G\in \mathcal{H}$, and an analogous equation holds for squared norms. 
\end{enumerate}}
\end{Remark}

This setting is suitable for introducing continuous analogs of sums of matricially free creation operators
$$
\wp(u)=\sum_{p,q=1}^r\wp_{p,q}(u),
$$
where the covariance of each $\wp_{p,q}(u)$ is assumed to be $b_{p,q}(u)\geq 0$, $u\in \mathpzc{U}$.
In particular, if $b_{p,q}(u)=d_{p}$ for any $p,q,u$, where $d_1+\cdots +d_r=1$, 
then $\{\wp(u):u\in \mathpzc{U}\}$ is a family of standard free creation operators. We would like to find a continuous analog of these decompositions, using direct integrals.

The continuous analogs of the matricially free creation operators will be denoted $\wp(f)$, where $f$ is an 
essentially bounded $\mathcal{G}$-valued function on $\Gamma_1$, 
namely $f\in L^{\infty}( \Gamma_{1}, \mathcal{G})$, where the square $\Gamma_1$ is 
equipped with the two-dimensional Lebesgue measure.

\begin{Definition}
{\rm For given $f\in L^{\infty}(\Gamma_{1}, \mathcal{G})$, let us define bounded linear 
operators $\wp(f)$ on $\mathcal{H}$ by 
$$
\wp(f)\left(\int^{\oplus}_{I}F_0(x_1)dx_1\right)=\int^{\oplus}_{\Gamma_{1}}f(x,x_1)F_0(x_1)dxdx_1
$$
for any $F_0\in L^{2}(I)$, and 
$$
\wp(f)\left(\int^{\oplus}_{\Gamma_{n}}F_n(x_1,\ldots ,x_{n+1})dx_1\ldots dx_{n+1}\right)
$$
$$
=
\int^{\oplus}_{\Gamma_{n+1}}f(x,x_1)\otimes F_{n}(x_1,\ldots , x_{n+1})dxdx_1\ldots dx_{n+1}
$$
for any $F_n\in L^{2}(\Gamma_{n}, {\mathcal G}^{\otimes n})$, where $n\in \mathbb{N}$. 
In the special case, when $f=g\otimes e(u)$, where $e(u)$ is some basis unit vector of ${\mathcal G}$, 
under the identification $L^{\infty}(\Gamma_1, {\mathcal G})\cong L^{\infty}(\Gamma_1)\otimes {\mathcal G}$, we will write 
$\wp(g,u)$ instead of $\wp(f)$.}
\end{Definition}

In order to give formulas for the adjoints of $\wp(f)$, we first need to define bounded operators 
which multiply each $F(\gamma)$ by an essentially bounded function $g$ of the first coordinate of $\gamma$. 
The explicit definition is given below.

\begin{Definition}
{\rm For $k\in L^{\infty}(I)$, define bounded linear operators
$$
M(k,\gamma):{\mathcal H}(\gamma)\rightarrow {\mathcal H}(\gamma)
$$ 
for any $\gamma=(x_1, \ldots , x_{n+1})\in \Gamma_{n}$ and $n\geq 0$ by
$$
M(k,\gamma)F_{n}(\gamma)=k(x_1)F_{n}(\gamma),
$$
and the associated decomposable operator in the direct integral form
$$
M(k):=\int_{\Gamma}^{\oplus}M(k,\gamma)d\gamma,
$$
which is a bounded linear operator on ${\mathcal H}$.
} 
\end{Definition}

The operator $M(k)$ reminds the gauge operator on the free Fock space 
associated with the multiplication operator by $k$, but one important difference is that 
$M(k)$ is non trivial on the vacuum space unless $k$ vanishes outside of the set of measure zero. 
Moreover, we will use the shorthand notations
\begin{eqnarray*}
F_{n}(\gamma)&=&f_{1}(x_1,x_2)\otimes \ldots \otimes f_{n}(x_n,x_{n+1}),\\
F_{n-1}(\gamma')&=&f_{2}(x_2,x_3)\otimes \ldots \otimes f_{n}(x_n,x_{n+1}),
\end{eqnarray*}
where $\gamma=(x_1, \ldots , x_{n+1})\in \Gamma_{n}$,
$\gamma'=(x_2, \ldots , x_{n+1})\in \Gamma_{n-1}$ and each $f(x_i,x_{i+1})$ is an element of
the Hilbert space $\mathcal{G}$.

\begin{Proposition}
The adjoints of the operators $\wp(f)$ are given by 
\begin{eqnarray*}
\wp^{*}(f)\int_{I}^{\oplus}F_{0}(\gamma)d\gamma_0&=&0\\
\wp^{*}(f)\int_{\Gamma_{n}}^{\oplus}F_{n}(\gamma)d\gamma_n&=
&\int_{\Gamma_{n-1}}^{\oplus} M(k,\gamma')F_{n-1}(\gamma')d\gamma_{n-1}'
\end{eqnarray*}
where 
$$
k(x_2)=\int_{0}^{1}\langle f_1(x_1,x_2),f(x_1,x_2)\rangle dx_{1},
$$ 
and $\langle., .\rangle $ is the canonical inner product in ${\mathcal G}$.
 
\end{Proposition}
{\it Proof.}
The first formula is obvious since the range of $\wp(f)$ is contained in the orthogonal complement 
of $L^{2}(I)$. To prove the second formula, we can take $F_{n}(\gamma)$ and $G_{n}(\gamma)$ 
to be simple tensors of the form
\begin{eqnarray*}
F_{n}(\gamma)&=&f_{1}(x_1,x_2)\otimes \ldots \otimes f_{n}(x_n,x_{n+1}),\\
G_{n}(\gamma)&=&g_{1}(x_1,x_2)\otimes \ldots \otimes g_{n}(x_n,x_{n+1}),
\end{eqnarray*}
where $\gamma=(x_1, \ldots , x_{n+1})$.
Then
\begin{eqnarray*}
&&
\langle \wp(f)\int_{\Gamma_{n-1}}^{\oplus}G_{n-1}(\gamma')d\gamma_{n-1}', 
\int_{\Gamma_{n}}^{\oplus}F_{n}(\gamma)d\gamma_{n}\rangle\\
&=&
\langle 
\int_{\Gamma_{n}}^{\oplus}f(x_1,x_2)\otimes G_{n-1}(\gamma')d\gamma_{n-1},
\int_{\Gamma_{n}}^{\oplus}F_{n}(\gamma)d\gamma_n\rangle\\
&=&
\int_{\Gamma_n}\langle f(x_1,x_2),f_1(x_1,x_2)\rangle \langle G_{n-1}(\gamma'),F_{n-1}(\gamma')\rangle dx_1\ldots dx_{n+1}\\
&=&
\int_{\Gamma_{n-1}}\left(\int_{I}\langle f(x_1,x_2),f_1(x_1,x_2)\rangle dx_1\right)\langle G_{n-1}(\gamma'), F_{n-1}(\gamma')\rangle dx_2\ldots dx_{n+1}\\
&=&
\langle \int_{\Gamma_{n-1}}^{\oplus}G_{n-1}(\gamma')d\gamma_{n-1}', \wp^{*}(f)\int_{\Gamma_{n}}^{\oplus}
F_{n}(\gamma)d\gamma_n(\gamma)\rangle,
\end{eqnarray*}
where $\gamma=(x_1, \ldots , x_{n+1})$ and $\gamma'=(x_2, \ldots , x_{n+1})$. The proof is completed.
\hfill $\blacksquare$

\begin{Corollary}
For any $f,f_1\in L^{\infty}(\Gamma_1, \mathcal{G})$, it holds that
$$
\wp^{*}(f)\wp(f_1)=M(k),
$$
where $k$ is of the same form as in Proposition 3.1.
\end{Corollary}

\begin{Remark}
{\rm Let us consider some special cases and one property of the operators studied above.
\begin{enumerate}
\item  
It is easy to see that if the functions $f,f_1$ do not depend on the second coordinate, i.e.
$f(x_1,x_2)=\widetilde{f}(x_1)$ and $f_1(x_1,x_2)=\widetilde{f_1}(x_1)$, then
$$
k(x_2)=\int_{0}^{1}\langle \widetilde{f_1}(x_1),\widetilde{f}(x_1)\rangle dx_{1},
$$
for any $x_2$ and thus $M(k)$ reduces to the multiplication by a constant and thus
we can write the relation
$$
\wp^{*}(f)\wp(f_1)=	\langle f_1,f\rangle=
\langle \widetilde{f_1},\widetilde{f}\rangle,
$$
and thus the operators $\wp(f), \wp^*(f)$ reduce to free creation and annihilation operators, respectively, 
with the natural inner product for square integrable ${\mathcal G}$-valued functions.
\item
In the above case, if we take two functions of the form: $f_1=\chi_{\Gamma_1}\otimes e(u')$ and $f_2=\chi_{\Gamma_1}\otimes e(u'')$, 
where $\chi_{\Gamma_1}$ is the characteristic function of the square,
and denote the associated creation operators by $\wp(u'),\wp(u'')$, respectively,  then it is easy to see that 
$\{\zeta(u):u\in \mathpzc{U}\}$, where
$$
\zeta(u)= \wp(u') +\wp^*(u'')
$$ 
and $u'\neq u''$, viewed as two `copies' of $u$, is a family of free circular operators 
(in other words, instead of the set $\mathpzc{U}$ we have to consider
a twice bigger set of indices).
\item
If we use characteristic functions of the triangle and take
$f=\chi_{\Delta}\otimes e(u)$ and $f_1=\chi_{\Delta}\otimes e(u_1)$, 
where $\Delta=\{(x,y): 0\leq x < y \leq 1\}$ and $e(u),e(u_1)$ are orthonormal basis vectors in ${\mathcal G}$, 
then
$$
k(x_2)=\delta_{u,u_1}\int_{0}^{x_2}dx_1=\delta_{u,u_1}x_2,
$$
and thus $M(k)$ reduces to the multiplication by $x_2$ times the Kronecker delta related to basis vectors, and thus the relation
between the creation and annihilation operators becomes
$$
\wp^{*}(f)\wp(f_1)=\delta_{u,u_1}M({\rm id}),
$$
which corresponds to the case when we deal with strictly upper triangular Gaussian random matrices and the operatorial 
limit is the triangular operator.
\item
For simplicity, we will assume from now on that $f=g\otimes e(u)$ and that $g$ does not depend on $u$.
Let us observe that if $(g_n)$ is a sequence of functions from $L^{\infty}(\Gamma_1)$ which converges
in norm to $g\in L^{\infty}(\Gamma_1)$, then the corresponding sequences of operators considered above 
converge strongly on ${\mathcal H}$, namely $s-\lim_{n\rightarrow \infty}\wp(g_n,u)=\wp(g,u)$,
$s-\lim_{n\rightarrow \infty}\wp^{*}(g_n,u)=\wp^{*}(g,u)$ and $s-\lim_{n\rightarrow \infty}M(g_n)=M(g)$.
\end{enumerate}}
\end{Remark}

\section{Continuous circular systems}

Other decompositions of ${\mathcal H}$ are also relevant since they 
give useful decompositions of the operators of interest.
We will introduce decompositions in which the sets of fibers are relatively small 
(indexed by $I$), but the fibers themselves are `long'.
These decompositions
allow us to introduce continuous analogs of matricial circular systems and 
interpret the operators of interest as integrals
of two-dimensional `densities'.

\begin{Definition}
{\rm For each $x\in [0,1]$, let us define the associated fiber Hilbert space that begins with $x$:
$$
\mathcal{N}(x):=\int_{\Gamma(y)}^{\oplus}{\mathcal N}(\gamma)d\widetilde{\gamma}
$$
$$
\cong {\mathbb C}\Omega(x)\oplus \int_{I}^{\oplus}{\mathcal H}(x,y)dy\oplus 
\int_{\Gamma_1}^{\oplus}{\mathcal H}(x,y)\otimes {\mathcal H}(y,z)dydz\oplus \ldots ,
$$
where $\Gamma(x)=\{(x, \gamma):\gamma\in \Gamma\}$ for any fixed $x\in I$, with 
$\widetilde{\gamma}(\{x\})=1$ and $\widetilde{\gamma}(\{x\}\times A)=\lambda(A)$ (the Lebesgue measure of $A$) 
for any $A\subset \Gamma$, and let
$$
{\mathcal H}=\int_{I}^{\oplus}{\mathcal N}(x)dx,
$$ 
be the associated direct integral decomposition. All Hilbert spaces involved are equipped with canonical 
inner products.}
\end{Definition}

\begin{Definition}
{\rm 
In a similar fashion, for all $(x,y)\in \Gamma_1$, define Hilbert spaces
$$
\mathcal{N}(x,y):={\mathcal H}(x,y)\oplus 
\int_{\Gamma_0}^{\oplus}{\mathcal H}(x,y)\otimes {\mathcal H}(y,z)dz\oplus \ldots ,
$$
equipped with the canonical inner products and let 
$$
{\mathcal H}\ominus {\mathcal H}_{0}=\int_{\Gamma_1}^{\oplus}{\mathcal N}(x,y)dxdy,
$$ 
be the associated direct integral decomposition, where Hilbert spaces involved are 
equipped with canonical inner products.
}
\end{Definition}

\begin{Definition}
{\rm Let us suppose that $\{e(y,z;u): u\in \mathpzc{U}\}$ 
is a countable orthonormal basis of ${\mathcal H}(y,z)$ for each 
$(y,z)\in \Gamma_{1}$. For any given $x,y\in I$ and $u\in \mathpzc{U}$, define isometries 
$\wp(x,y;u):{\mathcal N}(y)\rightarrow {\mathcal N}(x,y)$ 
by the direct integral extension of 
\begin{eqnarray*}
\wp(x,y;u)\Omega(y)&=&e(x,y;u),\\
\wp(x,y;u)e(y,z;s)&=&e(x,y;u)\otimes e(y,z;s),\\
\wp(x,y;u)(e(y,z;s)\otimes w)&=&e(x,y;u)\otimes e(y,z;s)\otimes w ,
\end{eqnarray*}
for any $x,y,z\in I$ and $u,s\in \mathpzc{U}$, where $e(y,z;s) \otimes w$ is a basis vector of 
some tensor product ${\mathcal H}(y,z)\otimes {\mathcal H}(z,z_1)\otimes \ldots\otimes {\mathcal H}(z_{n-1},z_{n})$.}
\end{Definition}

\begin{Remark}
{\rm 
Equivalently, we could act with $\wp(x,y;u)$ onto direct integrals in the last two equations. 
For instance, the second equation would then take the form
$$
\wp(x,y;u)\int_{I}^{\oplus}g(y,z)e(y,z;s)dz=\int_{I}^{\oplus}g(y,z)(e(x,y;u)\otimes e(y,z;s))dz,
$$
but it is more convenient to completely decompose the considered fibers since we get simpler formulas which are 
in correspondence with the discrete case. As far as this correspondence is concerned, in contrast to the 
discrete case, we do not include scalars in the definition of $\wp(x,y;u)$ in order to avoid lengthy formulas. 
These scalars, playing the role of covariances, are included in the function $g$ when we deal with 
$\wp(g,u)$ to the effect that $|g(x,y)|^2$ is the continuous analog of $b_{p,q}(u)$ (as we mentioned earlier, 
we shall assume for simplicity that these covariances do not depend on $u$).} 
\end{Remark}

\begin{Proposition}
For any $x,y\in I$ and $u\in \mathpzc{U}$, let $\wp^{*}(x,y;u):{\mathcal N}(x,y)\rightarrow {\mathcal N}(y)$ be the 
bounded operator defined by the direct integral extension of the formal formulas
\begin{eqnarray*}
\wp^{*}(x,y;u)e(x,y;u)&=&\Omega(y),\\
\wp^{*}(x,y;u)(e(x,y;u)\otimes w)&=&w ,
\end{eqnarray*}
for any $x,y\in I$ and $u\in \mathpzc{U}$ and $w$ as above, and setting them to be zero on the remaining basis vectors. 
Then the operator $\wp^{*}(x,y;u)$ is the adjoint of $\wp(x,y;u)$ for any $x,y,u$.
\end{Proposition}
{\it Proof.}
These formulas are obtained by straightforward computations.
\hfill $\blacksquare$

\begin{Definition}
{\rm Using the continuous family $\{\wp(x,y;u):x,y\in I, u\in \mathpzc{U}\}$ 
and the family of their adjoints, one then defines the continuous analogs of the matricial circular operators as
$$
\zeta(x,y;u)=\wp(x,y;u')+\wp^*(y,x,u''),
$$
for $(x,y)\in \Gamma_1$ and $u\in \mathpzc{U}$ and $u',u''$ are copies of $u$, as in Remark 3.2.
This definition is in agreement with the definition of matricial circular systems and therefore the family 
$$
\{\zeta(x,y;u): x,y\in I, \; u\in \mathpzc{U}\}
$$
will be called the {\it continuous circular system}. }
\end{Definition}

\begin{Proposition}
If $f(x,y)=g(x,y)\otimes e(u)$, the matrix elements of operators $\wp(f)=\wp(g,u)$ and their adjoints 
of the form
$$
\langle \wp(g,u)h_1,h_2\rangle=
\int_{\Gamma_1}g(x,y)\langle \wp(x,y;u)h_1(y), h_2(x,y)\rangle dxdy ,
$$
$$
\langle h_1,\wp^*(g,u)h_2\rangle=
\int_{\Gamma_1}\overline{g(x,y)}\langle h_1(y), \wp^*(x,y;u)h_2(x,y)\rangle dxdy ,
$$
where $g\in L^{\infty}(\Gamma_1)$ and $u\in \mathpzc{U}$, are well defined for any $h_1=\int_{I}^{\oplus}h_1(y)dy$ and $h_2=\int_{\Gamma_1}^{\oplus}h_2(x,y)dxdy$ 
according to the decompositions of ${\mathcal H}$ and ${\mathcal H}\ominus {\mathcal H}_{0}$, 
in Definitions 4.1 and 4.2, respectively. 
\end{Proposition}
{\it Proof.}
Observe that the integrals on the RHS are well defined since
$h_1$ and $h_2$ have square integrable norms by assumption, each $\wp(x,y;u)$ is an isometry from ${\mathcal N}(y)$ to
${\mathcal N}(x,y)$ and $g$ is essentially bounded on $\Gamma_1$. By Definition 3.2 and Proposition 3.1, 
the integrals on the RHS give the desired matrix elements. This completes the proof.
\hfill $\blacksquare$\\

In the above situation, we can write a decomposition of the creation 
operators $\wp(g,u)$ in the direct integral form
$$
\wp(g,u)=\int_{\Gamma_1}^{\oplus}g(x,y)\wp(x,y;u)dxdy ,
$$
and an analogous formula for the annihilation operators $\wp^*(g,u)$, namely
$$
\wp^*(g,u)=\int_{\Gamma_1}^{\oplus}\overline{g(x,y)}\wp^*(x,y;u)dxdy.
$$
We use the symbol $\oplus$ with a slight abuse of notation 
since the considered families of integrands 
are `almost decomposable' with respect to the direct integral decomposition of Definition 5.2.
The operators $\wp(x,y;u)$ ($\wp^*(x,y;u)$) can be interpreted as operators creating (annihilating) 
vector $e(u)$ of color $x$ `under condition $y$'.
The `condition' $y$ refers to the color of the vector onto which the operator $\wp(x,y;u)$ acts. 
If the given pairing is a block in the mixed *-moment of creation and annihilation operators, 
integration over $x$ of the associated $|g(x,y)|^2$ gives the contribution of the given pairing 
to the mixed *-moment.

It is also natural to define the corresponding 'circular operator' 
$$
\zeta(g,u)=\int_{\Gamma_1}^{\oplus}(g(x,y)\wp(x,y;u')+\overline{g(y,x)}\wp^*(y,x;u''))dxdy
$$
which becomes a circular operator $\zeta(u)$ if $g=\chi_{\Gamma_1}$.
Similarly, if we set $g=\chi_{\Gamma_1}$, we obtain canonical creation and annihilation operators associated with 
the basis vector $e(u)$, denoted $\wp(u)$ and $\wp^*(u)$.

\section{Mixed *-moments}

We would like to discuss the combinatorics of the mixed *-moments of the 
creation operators and of certain operators obtained from them.
A very interesting example is that of the triangular operator obtained 
as the limit realization of sctrictly upper triangular Gaussian random matrices. 

In our previous works we have studied the combinatorics of *-moments of 
various operators (creation, semicircular, circular, etc.) in the matricial (discrete) case [7,8,9,10].
It was based on the class of
{\it colored labeled noncrossing partitions} (if only one label is used, we just have colored noncrossing partitions). 
In the case of *-moments of creation, semicircular or circular operators, it suffices to consider pair partitions (see, for instance
[12]). The general idea of block coloring is very straightforward. We color the blocks with natural numbers from the finite 
set $[r]$, where $r$ is related to the block structure of the considered block random matrices in the sense that 
these matrices are assumed to have $r^2$ blocks. At the same time $r$ is the number of summands 
in the direct sum decomposition of our Fock space (with $r$ vacuum vectors). 
The color of each block is related in a matricial way to the color of its {\it nearest outer block} [7], as shown in Fig. ~1. 
In addition to colors, we equip blocks with labels, if necessary (namely, if we deal with more than one matrix), 
but labels are rather easy to deal with since they just have to match within a block. 

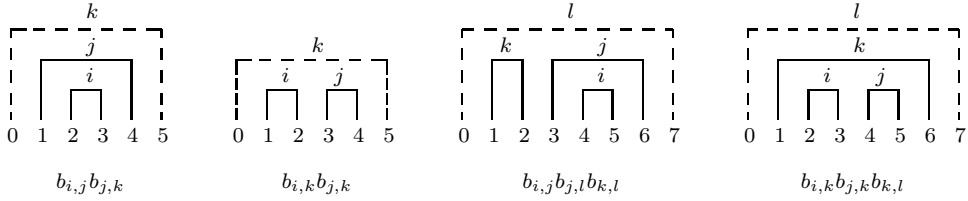
\begin{figure}
\unitlength=1mm
\special{em:linewidth 0.4pt}
\linethickness{0.4pt}
\begin{picture}(120.00,35.00)(30.00,0.00)
%%%%%%%%%%%%%diagram 1%%%%%%%%%%%%%%%%%%%%%%%%%%%%%%
\put(31.00,10.00){\line(0,1){8.00}}
\put(35.00,10.00){\line(0,1){4.00}}
\put(39.00,10.00){\line(0,1){4.00}}
\put(43.00,10.00){\line(0,1){8.00}}

\put(27.00,10.00){\line(0,1){1.00}}
\put(27.00,12.00){\line(0,1){1.50}}
\put(27.00,15.00){\line(0,1){1.50}}
\put(27.00,18.00){\line(0,1){1.50}}
\put(27.00,21.00){\line(0,1){1.00}}

\put(47.00,10.00){\line(0,1){1.00}}
\put(47.00,12.00){\line(0,1){1.50}}
\put(47.00,15.00){\line(0,1){1.50}}
\put(47.00,18.00){\line(0,1){1.50}}
\put(47.00,21.00){\line(0,1){1.00}}

\put(26.50,7.00){$\scriptstyle{0}$}
\put(30.50,7.00){$\scriptstyle{1}$}
\put(34.50,7.00){$\scriptstyle{2}$}
\put(38.50,7.00){$\scriptstyle{3}$}
\put(42.50,7.00){$\scriptstyle{4}$}
\put(46.50,7.00){$\scriptstyle{5}$}

\put(37.00,15.00){$\scriptstyle{i}$}
\put(37.00,19.00){$\scriptstyle{j}$}
\put(37.00,23.50){$\scriptstyle{k}$}

\put(27.00,22.00){\line(1,0){2.00}}
\put(30.00,22.00){\line(1,0){1.50}}
\put(33.00,22.00){\line(1,0){1.50}}
\put(36.00,22.00){\line(1,0){1.50}}
\put(39.00,22.00){\line(1,0){1.50}}
\put(42.00,22.00){\line(1,0){1.50}}
\put(45.00,22.00){\line(1,0){2.00}}

\put(31.00,18.00){\line(1,0){12.00}}
\put(35.00,14.00){\line(1,0){4.00}}
%%%%%%%%%%%%diagram 2%%%%%%%%%%%%%%%%%%%%%%%
\put(61.00,10.00){\line(0,1){4.00}}
\put(65.00,10.00){\line(0,1){4.00}}
\put(69.00,10.00){\line(0,1){4.00}}
\put(73.00,10.00){\line(0,1){4.00}}

\put(57.00,10.00){\line(0,1){1.50}}
\put(57.00,12.00){\line(0,1){1.50}}
\put(57.00,14.00){\line(0,1){1.50}}
\put(57.00,16.00){\line(0,1){2.00}}

\put(77.00,10.00){\line(0,1){1.50}}
\put(77.00,12.00){\line(0,1){1.50}}
\put(77.00,14.00){\line(0,1){1.50}}
\put(77.00,16.00){\line(0,1){2.00}}

\put(56.50,7.00){$\scriptstyle{0}$}
\put(60.50,7.00){$\scriptstyle{1}$}
\put(64.50,7.00){$\scriptstyle{2}$}
\put(68.50,7.00){$\scriptstyle{3}$}
\put(72.50,7.00){$\scriptstyle{4}$}
\put(76.50,7.00){$\scriptstyle{5}$}

\put(63.00,15.00){$\scriptstyle{i}$}
\put(70.00,15.00){$\scriptstyle{j}$}
\put(67.00,19.00){$\scriptstyle{k}$}

\put(57.00,18.00){\line(1,0){2.00}}
\put(60.00,18.00){\line(1,0){1.50}}
\put(63.00,18.00){\line(1,0){1.50}}
\put(66.00,18.00){\line(1,0){1.50}}
\put(69.00,18.00){\line(1,0){1.50}}
\put(72.00,18.00){\line(1,0){1.50}}
\put(75.00,18.00){\line(1,0){2.00}}

\put(61.00,14.00){\line(1,0){4.00}}
\put(69.00,14.00){\line(1,0){4.00}}
%%%%%%%%%%%%%diagram 3%%%%%%%%%%%%%%%%%%%%%%%%%%%%%%
\put(91.00,10.00){\line(0,1){8.00}}
\put(95.00,10.00){\line(0,1){8.00}}
\put(99.00,10.00){\line(0,1){8.00}}
\put(103.00,10.00){\line(0,1){4.00}}
\put(107.00,10.00){\line(0,1){4.00}}
\put(111.00,10.00){\line(0,1){8.00}}

\put(87.00,10.00){\line(0,1){1.00}}
\put(87.00,12.00){\line(0,1){1.50}}
\put(87.00,15.00){\line(0,1){1.50}}
\put(87.00,18.00){\line(0,1){1.50}}
\put(87.00,21.00){\line(0,1){1.00}}

\put(115.00,10.00){\line(0,1){1.00}}
\put(115.00,12.00){\line(0,1){1.50}}
\put(115.00,15.00){\line(0,1){1.50}}
\put(115.00,18.00){\line(0,1){1.50}}
\put(115.00,21.00){\line(0,1){1.00}}

\put(86.50,7.00){$\scriptstyle{0}$}
\put(90.50,7.00){$\scriptstyle{1}$}
\put(94.50,7.00){$\scriptstyle{2}$}
\put(98.50,7.00){$\scriptstyle{3}$}
\put(102.50,7.00){$\scriptstyle{4}$}
\put(106.50,7.00){$\scriptstyle{5}$}
\put(110.50,7.00){$\scriptstyle{6}$}
\put(114.50,7.00){$\scriptstyle{7}$}

\put(92.00,19.00){$\scriptstyle{k}$}
\put(105.00,15.00){$\scriptstyle{i}$}
\put(105.00,19.00){$\scriptstyle{j}$}
\put(101.00,23.50){$\scriptstyle{l}$}

\put(87.00,22.00){\line(1,0){1.50}}
\put(90.00,22.00){\line(1,0){1.50}}
\put(93.00,22.00){\line(1,0){1.50}}
\put(96.00,22.00){\line(1,0){1.50}}
\put(99.00,22.00){\line(1,0){1.50}}
\put(102.00,22.00){\line(1,0){1.50}}
\put(105.00,22.00){\line(1,0){1.50}}
\put(108.00,22.00){\line(1,0){1.50}}
\put(111.00,22.00){\line(1,0){1.50}}
\put(114.00,22.00){\line(1,0){1.00}}

\put(99.00,18.00){\line(1,0){12.00}}
\put(103.00,14.00){\line(1,0){4.00}}
\put(91.00,18.00){\line(1,0){4.00}}
%%%%%%%%%%%%diagram 4%%%%%%%%%%%%%%%%%%%%%%%
\put(129.00,10.00){\line(0,1){8.00}}
\put(133.00,10.00){\line(0,1){4.00}}
\put(137.00,10.00){\line(0,1){4.00}}
\put(141.00,10.00){\line(0,1){4.00}}
\put(145.00,10.00){\line(0,1){4.00}}
\put(149.00,10.00){\line(0,1){8.00}}

\put(125.00,10.00){\line(0,1){1.00}}
\put(125.00,12.00){\line(0,1){1.50}}
\put(125.00,15.00){\line(0,1){1.50}}
\put(125.00,18.00){\line(0,1){1.50}}
\put(125.00,21.00){\line(0,1){1.00}}

\put(153.00,10.00){\line(0,1){1.00}}
\put(153.00,12.00){\line(0,1){1.50}}
\put(153.00,15.00){\line(0,1){1.50}}
\put(153.00,18.00){\line(0,1){1.50}}
\put(153.00,21.00){\line(0,1){1.00}}

\put(124.50,7.00){$\scriptstyle{0}$}
\put(128.50,7.00){$\scriptstyle{1}$}
\put(132.50,7.00){$\scriptstyle{2}$}
\put(136.50,7.00){$\scriptstyle{3}$}
\put(140.50,7.00){$\scriptstyle{4}$}
\put(144.50,7.00){$\scriptstyle{5}$}
\put(148.50,7.00){$\scriptstyle{6}$}
\put(152.50,7.00){$\scriptstyle{7}$}

\put(135.00,15.00){$\scriptstyle{i}$}
\put(142.00,15.00){$\scriptstyle{j}$}
\put(139.00,19.00){$\scriptstyle{k}$}
\put(139.00,23.50){$\scriptstyle{l}$}

\put(125.00,22.00){\line(1,0){1.25}}
\put(127.67,22.00){\line(1,0){1.25}}
\put(130.33,22.00){\line(1,0){1.25}}
\put(132.98,22.00){\line(1,0){1.25}}
\put(135.63,22.00){\line(1,0){1.25}}
\put(138.30,22.00){\line(1,0){1.25}}
\put(140.95,22.00){\line(1,0){1.25}}
\put(143.60,22.00){\line(1,0){1.25}}
\put(146.25,22.00){\line(1,0){1.25}}
\put(149.00,22.00){\line(1,0){1.25}}
\put(151.70,22.00){\line(1,0){1.30}}

\put(129.00,18.00){\line(1,0){20.00}}
\put(133.00,14.00){\line(1,0){4.00}}
\put(141.00,14.00){\line(1,0){4.00}}
%%%%%%%%%%%%%%%%%%%%%%%%%%%%%%%%%%%%%%%
%%%%%%%%%%%%%%%%%%%%%%%%%%%%%%%%%%%%%%%
\put(33.00,1.00){$\scriptstyle{b_{i,j}b_{j,k}}$}
\put(63.00,1.00){$\scriptstyle{b_{i,k}b_{j,k}}$}
\put(95.00,1.00){$\scriptstyle{b_{i,j}b_{j,l}b_{k,l}}$}
\put(132.00,1.00){$\scriptstyle{b_{i,k}b_{j,k}b_{k,l}}$}

\end{picture}
\caption{Examples of colored non-crossing pair partitions. To each partition we assign the product
of covariances which depend on the colors of block and its nearest outer blocks.
Labels are assumed to be the same for all blocks and are omitted.}
\end{figure}

Roughly speaking, each creation-annihilation pair of operators indexed by the pair $(p,q)$ and labeled by $u$ 
produces the covariance $b_{p,q}(u)$, with $q$ depending on the vector in the space 
${\mathcal M}$ onto which the creation operator acts. If we associate a block with 
a given creation-annihilation pair, this vector corresponds to the nearest 
outer block of that block and must be colored by $q$ if the action is to be non-trivial.
That is why we color blocks with colors from the set $[r]$ and thus we can think that 
the covariance $b_{p,q}(u)$ is assigned to each block and its nearest outer block. 
Examples of such partitions are shown in Fig. 1, where, for simplicity, 
we assume that all labels are the same and can be omitted. One has to remark that we extend each partition by 
the imaginary block of some color $q$ which corresponds to the vacuum vector $\Omega_q$.
Now, the new idea in this paper is that in the limit $r\rightarrow \infty$ the combinatorics is 
still described by colored labeled noncrossing pair partitions, but the discrete set of colors $[r]$
is replaced by the interval $[0,1]$.

The basic definitions and notations are given below.
If $\pi$ is a non-crossing pair-partition of the set $[m]$, where $m$ is an even positive integer, 
which is denoted $\pi\in \mathcal{NC}_{m}^{2}$, the set 
$$
B(\pi)=\{V_1, \ldots , V_s\}
$$ 
is the set of its blocks, where $m=2s$. 
If $V_{i}=\{l(i),r(i)\}$ and $V_{j}=\{l(j),r(j)\}$ are two blocks of
$\pi$ with left legs $l(i)$ and $l(j)$ and right legs $r(i)$ and $r(j)$, respectively, then
$V_i$ is {\it inner} with respect to $V_j$ if $l(j)<l(i)<r(i)<r(j)$.
In that case $V_j$ is {\it outer} with respect to $V_i$.
It is the {\it nearest outer block} of $V_i$ if there is no block 
$V_k=\{l(k),r(k)\}$
such that $l(j)<l(k)<l(i)<r(i)<r(k)<r(j)$.
It is easy to see that the nearest outer block, if it exists, 
is unique, and we write in this case $V_j=o(V_i)$. If $V_i$ does not have an outer block, 
we set $o(V_i)=V_0$, where $V_{0}=\{0,m+1\}$ is the additional block called 
{\it imaginary}. The partition of the set $\{0,1, \ldots , m+1\}$ consisting of the blocks of $\pi$ and of the imaginary block 
will be denoted by $\widehat{\pi}$. 

Let us recall how we assigned colors to blocks in the discrete case [7,8,9,10], where we colored blocks with natural numbers from $[r]$. 
We used the set $\mathpzc{C}_r(\pi)$ of all mappings $c:B(\pi)\rightarrow [r]$ called {\it colorings}. 
By a {\it colored noncrossing pair partition} we understood a pair $(\pi,c)$, where $\pi \in \mathcal{NC}_{m}^{2}$ and 
$c\in \mathpzc{C}_r(\pi)$. Then, the set of pairs 
$$
B(\pi,f)=\{(V_1,c), \ldots , (V_s,c)\}
$$ 
played the role of the set of colored blocks. We assumed that also the imaginary block was colored by a number 
from the set $[r]$ and thus we could speak of a coloring of $\widehat{\pi}$. 
Then we assigned to blocks matricial elements associated with the covariance matrices
$B(u)=(b_{p,q}(u))\in M_{r}({\mathbb R})$, where $u\in \mathpzc{U}$.
For any $\pi\in \mathcal{NC}_{m}^{2}$ and $c\in \mathpzc{C}_{r}(\pi)$, let
$$
b_{q}(\pi, c)=\prod_{k=1}^{s}b_{q}(V_{k},c),
$$
where 
$$
b_{q}(V_k,c)=b_{s,t}(u)
$$
whenever $V_k=\{i,j\}$ is colored by $s$, its nearest outer block $o(V_k)$ is colored by $t$ and $u_i=u_j=u$ and
we assume that the imaginary block is colored by $q\in [r]$, and otherwise we set $b_{q}(V_k,c)=0$.
Examples of colored pair partitions with assigned weights are shown in Fig.~1. Let us just say that 
in the present paper we will repeat this procedure, except that the set of colors will be $x_1, \ldots , x_{s+1}$
lying in the interval $[0,1]$ and the weights will be products of $g(x_i,x_{o(i)})$, whenever $V_i$ is a block and $o(V_i)$ is 
its nearest outer block. 

In order to go from the discrete case to the continuous one, we will concentrate on the creation and annihilation 
operators. If a tuple $((\epsilon_1, u_1), \ldots, (\epsilon_m,u_m))$, where $\epsilon_j\in\{1,*\}$ and $u_j\in \mathpzc{U}$ for any $j$,
is given, where $m$ is even, we will say that $\pi\in \mathcal{NC}_{m}^{2}$ is {\it adapted} to it if $u_i=u_j$ 
whenever $\{i,j\}$ is a block and $(\epsilon_i, \epsilon_j)=(*,1)$ whenever $\{i,j\}$ is a block and $i<j$.
This notion (of a `noncrossing partition adapted to stars and labels')
is convenient when speaking of *-moments of free creation operators. Clearly,
if a tuple is given, it may have exactly one noncrossing pair partition adapted to it or 
none at all. If such a partition exists, then it means that the considered mixed *-moment of free 
creation operators gives a non-zero contribution.

We need to define a continuous analog of the state $\Psi$ considered in the discrete case. 
For simplicity, we can take $d_q=1/r$ for any $q$, which corresponds
to the decomposition of the random matrices into blocks which are asymptotically square and of equal sizes.  
We will use the state $\varphi:B({\mathcal H})\rightarrow {\mathbb C}$ of the form
$$
\varphi=\int_{I}^{\oplus}\varphi(\gamma)d\gamma,
$$
where $\varphi(\gamma)=\varphi(x)$ is the vacuum state 
associated with $\Omega(x)$, namely 
$$
\varphi(F)=\int_{I}\langle F(x)\Omega(x), \Omega(x)\rangle dx,
$$
where $F=\int_{I}^{\oplus}F(x)dx\in B({\mathcal H})$ according to the decomposition of ${\mathcal H}$ 
into fibers that end with $x\in I$, which is a natural continuous analog of the state $\Psi$ when $d_q=1/r$ 
for all $q$ obtained by taking the limit $r\rightarrow \infty$.

Computations of mixed *-moments of interest always reduce to the mixed *-moments of the 
{\it creation} operators. Therefore, let us first establish a connection on this 
level with the use of the operators $\wp(g,u)$ introduced in Section 3. 
We consider the case when $b_{p,q}(u)=d_p$ for any $p,q,u$. For further simplicity, 
one can even assume that $d_p=1/r$ for all $p$, but the result given below holds for any asymptotic 
dimensions.

\begin{Lemma}
For any $p_1,q_1, \ldots , p_m,q_m\in [r]$, $u_1, \ldots , u_m\in \mathpzc{U}$, 
$\epsilon_1, \ldots , \epsilon_m\in \{1,*\}$ and any $m\in \mathbb{N}$, it holds that
$$
\Psi(\wp_{p_1,q_1}^{\epsilon_1}(u_1)\ldots \wp_{p_m,q_m}^{\epsilon_m}(u_m))
=\varphi(\wp^{\epsilon_1}(f_1)\ldots \wp^{\epsilon_m}(f_m)),
$$   
where $f_k=g_k\otimes e(u_k)$ for $k\in [m]$ and $g_{k}$ is the characteristic function of the rectangle $I_{p_k}\times I_{q_k}\subset \Gamma_1$ for any $1\leq k\leq m$, where $I=I_1\cup\ldots \cup I_r$ is the partition of $I$ into disjoint non-empty intervals with natural ordering.
\end{Lemma}
{\it Proof.}
Let us observe that for any fixed $r\in \mathbb{N}$ an 
isometric embedding $\theta:{\mathcal M}\rightarrow {\mathcal H}$ is given by
\begin{eqnarray*}
\theta(\Omega_{q})&=&\frac{1}{\sqrt{d_q}}\int_{I_q}^{\oplus}\Omega(x)dx ,\\
\theta(e_{p_1,p_2}(u_1)\otimes \cdots \otimes e_{p_m,p_{n+1}}(u_n))&=&
\frac{1}{\sqrt{d_{p_1}\cdots d_{p_{n+1}}}}\int_{I_{p_1}\times \cdots \times I_{p_{n+1}}}^{\oplus}\\
&&
e(x_1,x_2,u_1)\otimes \cdots \otimes e(x_n,x_{n+1},u_n)dx_1\ldots dx_{n+1},
\end{eqnarray*}
for any $q,p_1,\ldots, p_{n+1}\in [r]$ and $u_1,\ldots , u_n\in \mathpzc{U}$.
It is then easy to check directly that the mixed *-moments of the operators 
$\wp_{p,q}(u)$ in the state $\Psi$ agree with the corresponding mixed *-moments 
of the operators $\wp(g,u)$, where $g$ is the characteristic function of $I_{p,q}$, 
respectively. This completes the proof. \hfill $\blacksquare$\\

The combinatorics of mixed *-moments of matricially free creation operators can be expressed in terms of 
noncrossing pair partitions adapted to stars and labels. 
It is not hard to see that they also describe the combinatorics of the mixed *-moments of 
the much more general family of operators $\wp(f)$, where $f=g\otimes e(u)$, in which 
matricially free creation operators are included if one takes characteristic functions of rectangles as above.
The main reason is that they encode two main facts: blocks must corresponds to pairings 
of creation and annihilation operators which have the same label, but the contribution of each partition depends
on the inner products.

\begin{Proposition}
Let $f_k=g_k\otimes e(u_k)$ and $\epsilon_k\in \{1,*\}$, 
where $k\in [m]$ and $m=2s$, be such that there exists a unique non-crossing 
pair partition $\pi\in \mathcal{NC}_{m}^{2}$ adapted to $((\epsilon_1,u_1) \ldots , (\epsilon_m,u_m))$.
Then
$$
\varphi(\wp^{\epsilon_1}(f_1)\ldots \wp^{\epsilon_m}(f_m))
=
\int_{\Gamma_{s}}\prod_{k=1}^{s}\langle f_{r(k)}(x_k, x_{o(k)}), f_{l(k)}(x_k,x_{o(k)})\rangle
dx_0dx_1\ldots dx_s ,
$$
where $V_k=\{l(k), r(k)\}$, $k=1, \ldots , s$, are the blocks of $\pi$ with $l(k)<r(k)$, with $x_{0}$ assigned 
to the imaginary block of $\pi$ and $\langle, .,.\rangle$ is the canonical inner product in ${\mathcal G}$.
\end{Proposition}
{\it Proof.}
Each pairing of a creation and annihilation operator
produces a function $g$ of one argument in the operator $M(k)$. Here, we just compute the inner products 
in ${\mathcal G}$ which appear in the definition of such $g$ for all pairings, which gives
$$
\langle f_{r(k)}(x_k, x_{o(k)}), f_{l(k)}(x_k,x_{o(k)})\rangle = 
\overline{g_{r(k)}(x_k,x_{o(k)})}g_{l(k)}(x_k, x_{o(k)})
$$ 
for each pairing, and then integrate the product of over all the variables $x_0, \ldots ,x_s$, which 
gives the desired formula.
\hfill $\blacksquare$

\begin{Proposition}
Under the above assumptions, if $g_{l(k)}=g_{r(k)}=\chi_{k}$ for all $k\in [s]$, where 
$\chi_1, \ldots , \chi_s$ are characteristic functions of some measurable 
subsets of $\Gamma_1$, then 
$$
\varphi(\wp^{\epsilon_1}(f_1)\ldots \wp^{\epsilon_m}(f_m))=Vol(\pi) ,
$$
where $Vol(\pi)$ is the volume of the region $V(\pi)\subseteq \Gamma_s$ defined by $x_0,x_1, \ldots , x_s$ for which
$\chi_k(x_k,x_{o(k)})=1$, for all $k=1, \ldots , s$, with $x_0$ assigned to the imaginary block. 
\end{Proposition}
{\it Proof.}
If we set the functions associated with the left and right legs of the block $\pi_k$ to be equal to
$\chi_k$, then in the proof of Proposition 5.1 we have
$$
\overline{g_{r(k)}(x_k,x_{o(k)})}g_{l(k)}(x_k, x_{o(k)}) 
=\chi_k(x_k,x_{o(k)}) ,
$$
which gives a condition on two variables, $x_k$ and $x_{o(k)}$, from among $s+1$ variables describing the $s+1$-dimensional cube 
$\Gamma_s$. Each inner product in the formula of Lemma 5.2 leads to a similar condition, which completes the proof.
\hfill  $\blacksquare$

\begin{Example}
{\rm
Consider the mixed *-moment associated with the pair partition $\pi=\{\{1,6\}, \{2,5\}, \{3,4\}\}$, where
it is natural to assume that $f_1=f_6$, $f_2=f_5$ and $f_3=f_4$. Then 
$$
\varphi(\wp^{*}(f_1)\wp^{*}(f_2)\wp^{*}(f_3)\wp(f_3)\wp(f_2)\wp(f_1))
$$
$$
=\int_{\Gamma_3}\parallel \!\!f_1(y,x)\!\!\parallel^2 \parallel \!\!f_2(z,y)\!\!\parallel^2 \parallel\!\! f_3(w,z)\!\!\parallel^2\,dwdzdydx.
$$
In the case when $f_j=g_j\otimes e(u_j)$ for $j=1,2,3$, we can replace the norms of $f_1,f_2,f_3$ by the absolute values of 
$g_1,g_2,g_3$, respectively. In particular, when these numerical valued functions
are the characteristic functions of the triangle $\Delta=\{(x,y): x\leq y\}$, this integral is equal 
to
$$
Vol(\pi)=\int_{0}^{1}dx\int_{0}^{x}dy\int_{0}^{y}dz\int_{0}^{z}dw=\frac{1}{24}.
$$
}
\end{Example}

\begin{Example}
{\rm
Consider the mixed *-moment associated with the pair partition $\pi=\{\{1,6\}, \{2,3\},\{4,5\}\}$, where 
it is natural to assume that $f_1=f_6$, $f_2=f_3$ and $f_4=f_5$. Then 
$$
\varphi(\wp^{*}(f_1)\wp^*(f_2)\wp(f_2)\wp^*(f_4)\wp(f_4)\wp(f_1))
$$
$$
=
\int_{\Gamma_3}\parallel \!\!f_1(y,x)\!\!\parallel^2 \parallel\!\! f_2(z,y)\!\!\parallel^2
\parallel \!\!f_4(w,y)\!\!\parallel^2
\,dwdzdydx.
$$
In the case when $f_j=g_j\otimes e(u_j)$ for $j=1,2,4$ and the numerical valued functions are the characteristic 
functions of the same triangle $\Delta$ as in Example 5.1, this integral is equal to
$$
Vol(\pi)=\int_{0}^{1}dx\int_{0}^{x}dy\int_{0}^{y}dz\int_{0}^{y}dw=\frac{1}{12}.
$$
It is not a coincidence that this volume is twice bigger than that in Example 5.1 
(in the former case we had one $4$-dimensional simplex, here we have the union of two such 
simplices).}
\end{Example}

Finally, let us return to the asymptotic *-distributions of Gaussian random matrices
with i.b.i.d. entries. We assume that we have $r^2$ blocks for each $r\in {\mathbb N}$.
Later we will go with $r$ to infinity. Therefore, at this point it seems appropriate to 
include $r$ in the symbols denoting random matrices as well as limit operators.

\begin{Proposition}
For any $r\in {\mathbb N}$, let $\{Y(u,n,r): u\in \mathpzc{U}\}$ be a family of square 
$n\times n$ independent complex Gaussian random matrices with i.b.i.d. entries for any natural $n$ 
and any $u\in \mathpzc{U}$. Then,    
$$
\lim_{n\rightarrow \infty}
\tau(n)(Y^{\epsilon_1}(u_1,n,r)\ldots Y^{\epsilon_m}(u_m,n,r))=
\Psi(\zeta^{\epsilon_1}(u_1,r)\ldots \zeta^{\epsilon_m}(u_m,r))
$$
for any $u_1, \ldots , u_m\in \mathpzc{U}$ and $\epsilon(1), \ldots , \epsilon(m)\in \{1,*\}$, where
$$
\zeta(u,r)=\sum_{p,q=1}^{r}\zeta_{p,q}(u,r)
$$
for any $u\in \mathpzc{U}$ and the operators $\zeta_{p,q}(u,r)$ are
discrete (generalized) matricial circular operators corresponding to given $r$.
\end{Proposition}
{\it Proof.}
The proof of this result was given in [10].
\hfill $\blacksquare$\\

The next step consists in taking the limit of the *-moments on the RHS as $r\rightarrow \infty$.
We assume that $d_{p}=1/r$ for any $p$ and any $r$ and, for simplicity, we assume that 
the block covariances $b_{p,q}(u,r)$, defined for all $p,q\in [r]$ and 
all $u,r$, do not depend on $u$. From these block covariances we built a sequence of 
simple functions
$$
b_r(x,y)=\sum_{p,q=1}^{r}b_{p,q}(u,r)\chi_{I_p,I_q}(x,y)
$$
and assume that it converges to some $g\in L^{\infty}(\Gamma_1)$ 
as $r\rightarrow \infty$. Then we compute the limit of the mixed *-moments expressed as liner combinations of
*-moments of the type given by Lemma 5.1. 

\begin{Theorem}
Let $\{Y(u,n,r):u\in \mathpzc{U}, r\in {\mathbb N}\}$ be a family of independent $n\times n$ 
random matrices for any $n\in \mathbb{N}$, such that 
\begin{enumerate}
\item
each $Y(u,n,r)$ consists of $r^2$ blocks of equal size with i.b.i.d. complex Gaussian entries,
\item
the sequence of simple functions $(b_r)$ converges to $g$ 
in $L^{\infty}(\Gamma_1)$ as $r\rightarrow \infty$.
\end{enumerate}
Then 
$$
\lim_{r\rightarrow \infty}\lim_{n\rightarrow \infty}
\tau(n)(Y^{\epsilon_1}(u_1,n,r)\ldots Y^{\epsilon_m}(u_m,n,r))=
\varphi(\eta^{\epsilon_1}(g,u_1)\ldots \eta^{\epsilon_m}(g,u_m)) ,
$$
where $\eta(g,u_j)=\wp(g,u'_j)+\wp^*(g^t,u''_j)$, with $g^t(x,y)=g(y,x)$ and all labels 
$u'_j, u''_j$, $j\in [m]$, different, and where $\varphi=\int_{I}^{\oplus}\varphi(\gamma)d\gamma$.
\end{Theorem}
{\it Proof.}
The second limit ($n\rightarrow \infty$) was computed in Proposition 5.4. 
It is easy to see that the moments obtained there, namely $\Psi(\zeta^{\epsilon_1}(u_1,r)\ldots \zeta^{\epsilon_m}(u_m,r))$,
can be written as linear combinations of mixed *-moments of creation and annihilation operators 
of continuous type in the state $\varphi$, namely such as those given in Lemma 5.1, since 
$$
\zeta(u,r)=\sum_{p,q=1}^{r}\zeta_{p,q}(u,r)=\sum_{p,q}^{r}(\wp_{p,q}(u',r)+\wp_{q,p}^*(u'',r)) ,
$$
where the matricial creation and annihilation operators are assumed to have covariances independent of $u$, but
otherwise arbitrary nonnegative numbers, i.e. $b_{p,q}(u)=b_{p,q}$ for any $u$. When we express the RHS 
in terms of operators of the form $\wp(g_{p,q},u)$ and their adjoints, where $g_{p,q}=b_{p,q}\chi_{I_p\times I_q}$ for any $p,q$,
we can write the above sum as  
$$
\eta(g_r,u)=\wp(g_r,u')+\wp^*(g^t_r,u'') ,
$$
where $g_r(x,y)=\sum_{p,q=1}^{r}g_{p,q}$ and $g^t_r$ stands for the transpose of $g_r$. 
Now, if $(g_r)$ converges to $g$ in $L^{\infty}(\Gamma_1)$,
then, by Definition 3.2, the mixed *-moments of $\wp(g_r,u)$ in the state $\varphi$ converge to the corresponding mixed *-moments of 
$\wp(g,u)$, which entails convergence of the mixed *-moments of $\eta(g_r,u)$ in the state $\varphi$ 
to the mixed *-moments of $\eta(g,u)$. This completes the proof. 
\hfill $\blacksquare$

\section{Triangular operator and labeled ordered trees}

Let us apply Theorem 5.1 to independent strictly upper triangular Gaussian random matrices, whose limits are free
triangular operators [4]. We express the limit *-moments of such matrices in terms of operators $\eta(u)=\zeta(\chi_{\Delta},u)$
where $\chi_{\Delta}$ is the characteristic function of the triangle $\Delta$. Our result gives a new Hilbert space realization 
of the limit *-moments, equivalent to the von Neumann algebra approach of Dykema and Haagerup in [4], 
where the triangular operator $T$ was introduced. Note that our approach to the combinatorics of its *-moments
is also different than the algorithm in [4, Lemma 2.4]. 

By a family of independent strictly upper triangular Gaussian random matrices we will understand a 
family of complex $n\times n$ matrices $Y(u,n)$, where $u\in \mathpzc{U}$, whose entries above the main diagonal 
form a family of complex Gaussian random variables whose real and imaginary parts form a family of 
$n(n-1)$ i.i.d. Gaussian random variables for each $u$ (also independent for different $u\in \mathpzc{U}$), 
each having mean zero and variance $1/2n$. 

\begin{Theorem}
Let $\{Y(u,n):u\in \mathpzc{U}\}$ be a family of independent strictly  
upper triangular Gaussian random matrices for any $n\in \mathbb{N}$. Then
$$
\lim_{n\rightarrow \infty}
\tau(n)(Y^{\epsilon_1}(u_1,n)\ldots Y^{\epsilon_m}(u_m,n))=
\varphi(\eta^{\epsilon_1}(u_1)\ldots \eta^{\epsilon_m}(u_m))
$$
for any $\epsilon_1, \ldots, \epsilon_m\in \{1,*\}$ and $u_1, \ldots , u_m\in \mathpzc{U}\}$, 
where $\eta(u_j)=\eta(\chi_{\Delta},u_j)$, $j\in [m]$, with  
$\Delta=\{(x,y)\in \Gamma_1: x<y\}$ and $\varphi=\int_{I}^{\oplus}\varphi(\gamma)d\gamma$.
\end{Theorem}
{\it Proof.}
We know from [4] that the limits of the mixed *-moments of independent strictly upper triangular Gaussian 
random matrices under $\tau(n)$ as $n\rightarrow \infty$ exist and are, by definition, 
equal to the mixed *-moments of free triangular operators $T(u)$, namely
$$
\lim_{n\rightarrow \infty}\tau(n)(Y^{\epsilon_1}(u_1,n)\cdots Y^{\epsilon_m}(u_m,n))=
\varphi(T^{\epsilon_1}(u_1)\cdots T^{\epsilon_m}(u_m))
$$
for any $\epsilon_1, \ldots , \epsilon_m\in \{1,*\}$ and any $u_1, \ldots , u_m\in \mathpzc{U}\}$.
At the same time, it can be justified that the above limits are equal 
to the limit moments of Theorem 5.1, where matrices $\{Y(u,n,r): u\in \mathpzc{U}\}$ 
are $n\times n$ independent block strictly upper triangular matrices with $r^2$ blocks for all natural $n$ and $r$, 
the non-vanishing blocks being $S_{p,q}(u,n,r)$ for $p<q$. For instance, an estimate in terms of
Schatten $p$-norms $\|A\|_{p}=\sqrt[p]{{\rm tr}(n)(|A|^{p})}$ for $p\geq 1$, where ${\rm tr}(n)$ is the 
normalized trace, can be used. It holds that 
$$
|{\rm tr}(n)(A)|\leq \parallel \!A \!\parallel _1\leq \parallel \! A \!\parallel _{p}\leq \parallel \!A\!\parallel
$$
for any $p\geq 1$. Therefore, let $Y_j=Y^{\epsilon_j}(u_j,n)$ and $Y_j'=Y^{\epsilon_j}(u_j,n,r)$ for $j\in [m]$ and any $n,r$. Then,
applying the above inequalities to the trace
$$
{\rm tr}(n)(Y_1\ldots Y_m-
Y_1'\ldots Y_m')
=\sum_{j=1}^{m}
{\rm tr}(n)(Y_1'\ldots Y_{j-1}'(Y_j-Y_j')Y_{j+1}\ldots Y_m) ,
$$
and using repeatedly the H\"{o}lder inequality 
$\parallel \!AB \!\parallel_{s}\leq \parallel \!A \!\parallel_{p}\parallel \!B \!\parallel_{q}$, where 
$s^{-1}=p^{-1}+q^{-1}$, we obtain an upper bound for the absolute value of this trace 
of the form
$$
(m+1) \cdot \max_{1\leq j \leq m}\|Y_j-Y_j'\|_{m+1}\cdot (\max_{1\leq k \leq m}\{\|Y_k\|_{m}, \|Y_{k}'\|_{m}\})^{m-1}.
$$
Therefore, for large $n$ there exists $R$ such that if $r>R$, then the difference between the mixed 
*-moments of primed and unprimed matrices can be made arbitrarily small
since the norms $\|Y_j-Y_j'\|_{m+1}$ can be made arbitrarily small for large $n$ and large $r>R$.
Therefore, we can use Theorem 5.1 to express the limit mixed *-moments of the strictly upper 
triangular matrices in terms of the operators $\eta(\xi_{\Delta}, u)$, respectively, where 
$\xi_{\Delta}$ is the characteristic function of the triangle $\Delta=\{(x,y)\in \Gamma_1: x<y\}$.
In other words, we can identify the triangular operators $T(u)$ with 
$\eta(\xi_{\Delta}, u)$, $u\in \mathpzc{U}$. This completes the proof.
\hfill $\blacksquare$\\

In order to give these *-moments in a more explicit form, let us assign continuous colors to 
blocks of $\pi\in \mathcal{NC}^2((\epsilon_1, u_1), \ldots , (\epsilon_m,u_m))$, where $m=2s$, and 
analyze relation between these colors. By $\mathcal{NC}^2((\epsilon_1, u_1), \ldots , (\epsilon_m,u_m))$ 
we denote the set of noncrossing pair partitions of $[m]$, such that $u_i=u_j$ and 
$\epsilon_i \neq \epsilon_j$ whenever $\{i,j\}$ is a block.
We should remember that stars refer here to operators of the form
$$
\eta(f,u)=\wp(f,u')+\wp^{*}(f^t,u''),
$$
where $f=\chi_{\Delta}$. For simplicity, let us write $\eta(f,u)=\eta(u)$, $\wp(f,u')=\wp(u')$ 
and $\wp(f^t,u'')=\wp(u'')$. 
Each pairing that gives a nonzero contribution must be of the form $(\wp^*(u'),\wp(u'))$ or 
$(\wp^*(u''), \wp(u''))$. 
The first one is obtained when $\eta^*(u)$ is associated with the left leg of a block and $\eta(u)$ is associated 
with the right leg, whereas in the second one the stars are interchanged. In any case, 
only one leg of a block can be marked with a star. We do not star the legs of the imaginary block.

\begin{Definition}
{\rm 
Let $V_j$ be a block of $\pi\in \mathcal{NC}^2((\epsilon_1, u_1), \ldots , (\epsilon_m,u_m))$ 
and let $V_{o(j)}$ be its nearest outer block. We distinguish four types of blocks:
\begin{enumerate}
\item {\it type 1}: the right leg of $V_j$ and the left leg of $V_{o(j)}$ are starred, 
\item {\it type 2}: the left leg of $V_j$ and the right leg of $V_{o(j)}$ are starred,
\item {\it type 3}: the left legs of both $V_j$ and $V_{o(j)}$ are starred, 
\item {\it type 4}: the right legs of both $V_j$ and $V_{o(j)}$ are starred.
\end{enumerate}
All types of pairs $(V,o(V))$ are shown in Fig.~2.}
\end{Definition}

\begin{Remark}
{\rm Our combinatorics is based on coloring the blocks of noncrossing pair partitions with numbers 
from $[0,1]$ and finding relations betwen them. This is a continuous analog of coloring blocks
with numbers from the discrete set $[r]$.
\begin{enumerate}
\item
Let us color the blocks of $\pi\in \mathcal{NC}^2((\epsilon_1, u_1), \ldots , (\epsilon_m,u_m))$, where $m=2s$,
and the imaginary block with $s+1$ continuous colors form $[0,1]$: 
$x_1, \ldots, x_{s+1}$, assigned to $V_1, \ldots , V_{s+1}$, respectively. 
It is convenient to number those blocks and colors starting from the right, as shown in Fig.~4 (thus, for instance,
$V_1$ is the imaginary block and its color is $x_1$). 
\item
Now, using these inequalities, we can associate a region $V(\pi)\subset \Gamma_2$ to each $\pi$. Namely, 
let 
$$
V(\pi)=\{x\in\Gamma_s:\; x_j<x_{o(j)}\; {\rm if}\;V_j\in B'(\pi)\;\wedge\; x_j>x_{o(j)}\; {\rm if}\;
V_j\in B''(\pi)\},
$$
where $B'(\pi)$ and $B''(\pi)$ stand for blocks of $\pi$ whose left legs are starred and unstarred, respectively.
\end{enumerate}
}
\end{Remark}

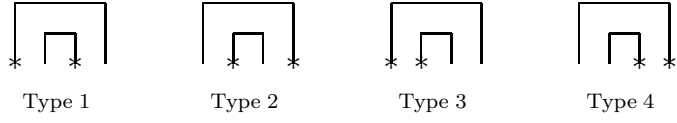
\begin{figure}
\unitlength=1mm
\special{em.linewidth 0.5pt}
\linethickness{0.5pt}
\begin{picture}(80.00,30.00)(0.00,-08.00)

%%%%%%%%partition 1%%%%%%%%%%%%%%%%
\put(00.00,00.00){\line(0,1){8.00}}
\put(04.00,00.00){\line(0,1){4.00}}
\put(08.00,00.00){\line(0,1){4.00}}
\put(12.00,00.00){\line(0,1){8.00}}
\put(0.00,08.00){\line(1,0){12.00}}
\put(04.00,04.00){\line(1,0){4.00}}

\put(-01.00,-02.00){*}
\put(07.00,-02.00){*}
\put(01.00,-06.00){$\scriptstyle{{\rm Type\;1}}$}
%%%%%%%%partition 2%%%%%%%%%%%%%%%%
\put(25.00,00.00){\line(0,1){8.00}}
\put(29.00,00.00){\line(0,1){4.00}}
\put(33.00,00.00){\line(0,1){4.00}}
\put(37.00,00.00){\line(0,1){8.00}}
\put(25.00,08.00){\line(1,0){12.00}}
\put(29.00,04.00){\line(1,0){4.00}}

\put(28.00,-02.00){*}
\put(36.00,-02.00){*}
\put(26.00,-06.00){$\scriptstyle{{\rm Type\;2}}$}
%%%%%%%%partition 3%%%%%%%%%%%%%%%%
\put(50.00,00.00){\line(0,1){8.00}}
\put(54.00,00.00){\line(0,1){4.00}}
\put(58.00,00.00){\line(0,1){4.00}}
\put(62.00,00.00){\line(0,1){8.00}}
\put(50.00,08.00){\line(1,0){12.00}}
\put(54.00,04.00){\line(1,0){4.00}}

\put(49.00,-02.00){*}
\put(53.00,-02.00){*}
\put(51.00,-06.00){$\scriptstyle{{\rm Type\;3}}$}
%%%%%%%%partition 4%%%%%%%%%%%%%%%%
\put(75.00,00.00){\line(0,1){8.00}}
\put(79.00,00.00){\line(0,1){4.00}}
\put(83.00,00.00){\line(0,1){4.00}}
\put(87.00,00.00){\line(0,1){8.00}}
\put(75.00,08.00){\line(1,0){12.00}}
\put(79.00,04.00){\line(1,0){4.00}}

\put(82.00,-02.00){*}
\put(86.00,-02.00){*}
\put(76.00,-06.00){$\scriptstyle{{\rm Type\;4}}$}

\end{picture}
\caption{Four types of pairs $(V, o(V))$, where $V$ is a block and $o(V)$ is the nearest outer block of $V$, which 
depend on which legs are starred.} 
\end{figure}

\begin{Corollary}
The non-vanishing mixed *-moments of the free triangular operators in the state $\varphi$ take the form 
$$
\varphi(T^{\epsilon_1}(u_1)\ldots T^{\epsilon_m}(u_m))=
\sum_{\pi\in \mathcal{NC}^{2}((\epsilon_1, u_1), \ldots, (\epsilon_m, u_m))}
Vol(\pi),
$$
where $m=2s$ and $Vol(\pi)$ is the volume of the region $V(\pi)$.
\end{Corollary}
{\it Proof.}
We know that we can replace $T(u_j)$ by $\eta(u_j)$. Now, without loss of generality, 
we can assume that $u_j=u$ for all $j$ since the general case just gives the additional condition on $\pi$ 
that $u_i=u_j$ whenever $\{i,j\}$ is a block. We will omit $u$ and write 
$\wp_1=\wp(u')$ and $\wp_2=\wp(u'')$ and replace $\mathcal{NC}^{2}((\epsilon_1, u), \ldots, (\epsilon_m, u))$ by 
$\mathcal{NC}^{2}(\epsilon_1, \ldots, \epsilon_m)$.
We have
$$
\varphi(T^{\epsilon_1}\ldots T^{\epsilon_m})=
\sum_{\pi\in \mathcal{NC}^{2}(\epsilon_1, \ldots, \epsilon_m)}
\varphi(\wp_{j_1(\pi)}^{\epsilon_1(\pi)}\ldots \wp_{j_m(\pi)}^{\epsilon_m(\pi)}),
$$
where the partition $\pi$ `chooses' whether to take the pairing $(\wp_1^*,\wp_1)$ or 
$(\wp_2^*,\wp_2)$, which formally can be written as
$$
\wp_{j_i(\pi)}^{\epsilon_i(\pi)}=\left\{
\begin{array}{ll}
\wp_1& {\rm if}\;\epsilon_i=1\;{\rm and}\; i\in \mathcal{R}(\pi)\\
\wp_2& {\rm if}\;\epsilon_i=*\;{\rm and}\; i\in \mathcal{R}(\pi)\\
\wp_2^*& {\rm if}\;\epsilon_i=1\;{\rm and}\; i\in \mathcal{L}(\pi)\\
\wp_1^*& {\rm if}\;\epsilon_i=*\;{\rm and}\; i\in \mathcal{L}(\pi)
\end{array}
\right.,
$$
where $\mathcal{R}(\pi)$ and $\mathcal{L}(\pi)$ stand for the right and left legs of $\pi$, respectively.
If $j>1$, then we choose the color $x_j$ assigned to block $V_j$ to be the first coordinate of $\chi_{\Delta}$ associated with 
the pairing of type $(\wp_1^*, \wp_1)$, or the first coordinate of $\chi_{\Delta}^{t}$ associated with the pairing of 
type $(\wp_2^*, \wp_2)$, depending on whether the left leg of $V_j$ is starred or unstarred, respectively.
Since we have $f=\chi_{\Delta}$ in each operator $\eta$ that appears in the moment
$\varphi(\eta^{\epsilon_1}\ldots \eta^{\epsilon_m})$, let us observe that if 
the left leg of $V_j$ is starred, then $x_j<x_{o(j)}$ is obtained from Remark 3.3 on the form 
of $M(k)$, with $k(x_{o(j)})=\int_{x_j<x_{o(j)}}dx_j$. 
In turn, if the left leg of $V_j$ is unstarred, then instead of $\chi_{\Delta}$, we 
take its transpose in the corresponding pairing of type $(\wp_2^*, \wp_2)$
which amounts to taking $M(k)$ with $k(x_{o(j)})=\int_{x_j>x_{o(j)}}dx_j$ (since
$u_j=u_{o(j)}$, by the adaptedness assumption on $\pi$). This completes the proof.
\hfill $\blacksquare$\\

\begin{figure}
\unitlength=1mm
\special{em.linewidth 0.5pt}
\linethickness{0.5pt}
\begin{picture}(140.00,35.00)(-20.00,-5.00)
%tree 1%%%%%%%%%%%%%%%%%%%%%%%%%%%%%
\put(10.00,05.00){\line(0,1){15.00}}
\put(10.00,5.00){\circle*{1.00}}
\put(10.00,10.00){\circle*{1.00}}
\put(10.00,15.00){\circle*{1.00}}
\put(10.00,20.00){\circle*{1.75}}
\put(9.00,0.00){$\scriptstyle{T_{1}}$}
%tree 2%%%%%%%%%%%%%%%%%%%%%%%%%%%%%
\put(25.00,15.00){\line(0,1){5.00}}
\put(25.00,15.00){\line(1,-1){5.00}}
\put(25.00,15.00){\line(-1,-1){5.00}}
\put(30.00,10.00){\circle*{1.00}}
\put(20.00,10.00){\circle*{1.00}}
\put(25.00,15.00){\circle*{1.00}}
\put(25.00,20.00){\circle*{1.75}}
\put(24.00,0.00){$\scriptstyle{T_{2}}$}
%tree 3%%%%%%%%%%%%%%%%%%%%%%%%%%%%%
\put(45.00,20.00){\line(1,-1){5.00}}
\put(45.00,20.00){\line(-1,-1){5.00}}
\put(40.00,15.00){\line(0,-1){5.00}}
\put(40.00,15.00){\circle*{1.00}}
\put(40.00,10.00){\circle*{1.00}}
\put(50.00,15.00){\circle*{1.00}}
\put(45.00,20.00){\circle*{1.75}}
\put(44.00,0.00){$\scriptstyle{T_{3}}$}
%tree 4%%%%%%%%%%%%%%%%%%%%%%%%%%%%%
\put(65.00,20.00){\line(1,-1){5.00}}
\put(65.00,20.00){\line(-1,-1){5.00}}
\put(70.00,15.00){\line(0,-1){5.00}}
\put(70.00,10.00){\circle*{1.00}}
\put(60.00,15.00){\circle*{1.00}}
\put(70.00,15.00){\circle*{1.00}}
\put(65.00,20.00){\circle*{1.75}}
\put(64.00,0.00){$\scriptstyle{T_{4}}$}
%tree 5%%%%%%%%%%%%%%%%%%%%%%%%%%%%%
\put(85.00,20.00){\line(0,-1){5.00}}
\put(85.00,20.00){\line(1,-1){5.00}}
\put(85.00,20.00){\line(-1,-1){5.00}}
\put(80.00,15.00){\circle*{1.00}}
\put(85.00,15.00){\circle*{1.00}}
\put(90.00,15.00){\circle*{1.00}}
\put(85.00,20.00){\circle*{1.75}}
\put(84.00,0.00){$\scriptstyle{T_{5}}$}
\end{picture}
\caption{Ordered rooted trees on 4 vertices}
\end{figure}
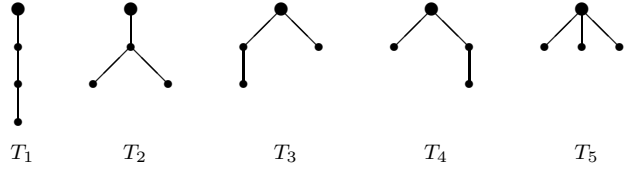

\begin{Remark}
{\rm Let us make some remarks on noncrossing pair partitions, colored noncrossing pair partitions, 
ordered rooted trees and labeled ordered rooted trees, which will be useful in establishing 
a bijective proof of the formula for the moments $\varphi((T^*T)^n)$ for the triangular operator $T$.
\begin{enumerate}
\item
There is an obvious bijection between $\mathcal{NC}_{2n}^{2}$ and
the set of associated extended pair-partitions of $[0,2n+1]$, and for that reason 
we will use the same notation for this set and we will understand from now on that each
$\pi\in {\mathcal NC}_{2n}^{2}$ is identified with its extension $\pi\cup \{0,2n+1\}$.
\item
For any integer $n$, let $\mathcal{O}_{n}$ be the set of ordered rooted trees on the set of $n+1$ vertices (see Example 6.1).
There is a natural bijection
$$
\gamma:\mathcal{O}_{n}\rightarrow \mathcal{NC}_{2n}^{2}
$$
given by the following rule: a vertex $v$ of $\mathcal{O}_{n}$ is a child of vertex $w$ 
if and only if $\gamma(w)$ is the nearest outer block of $\gamma(v)$. 
Thus, the root of $T$ corresponds to the imaginary block of $\pi=\gamma(T)$. 
\item
Suppose now that to each vertex of an ordered rooted tree we assign a {\it label} from some set 
$\mathcal{L}$. We will consider the case when trees on $n+1$ vertices are labeled 
by numbers from the set $[n+1]$. The bijection between noncrossing pair partitions and 
ordered rooted trees leads to a natural bijection 
$$
\gamma':\mathcal{L}_{n}\rightarrow \mathcal{CNC}_{2n}^{2}
$$
between labeled ordered rooted trees on $n+1$ vertices, denoted ${\mathcal L}_{n}$, 
where the vertices are labeled by different numbers from the set $[n+1]$, and
colored non-crossing pair-partitions of the set $[0,2n+1]$, denoted 
$\mathcal{CNC}_{2n}^{2}$, extended by the (colored) imaginary block (again, we identify $\pi$ with $\pi \cup\{0,2n+1\}$), where blocks (including the imaginary block) are colored by different numbers from the set $[n+1]$. 
Clearly, the labeling of each vertex is the same as the coloring of the corresponding block. We prefer 
to speak of {\it labeled} ordered trees and {\it colored} noncrossing partitions since the first terminology is 
standard in the combinatorics of trees and we used the second one in our previous works.
\item
It is well known that the number of ordered labeled trees on $n+1$ vertices 
is given by 
$$
|\mathcal{L}_{n}|=(n+1)!C_n=\frac{(2n)!}{n!},
$$
where $C_n=\frac{1}{n+1}{2n \choose n}$ is the $n$th Catalan number. 
\end{enumerate}
}
\end{Remark}

\begin{figure}
\unitlength=1mm
\special{em.linewidth 0.5pt}
\linethickness{0.5pt}
\begin{picture}(140.00,155.00)(-02.00,10.00)

\put(24.00,156.00){$\pi$}
\put(55.00,156.00){$R(\pi)$}
\put(90.00,156.00){$S(\pi)$}
\put(126.00,156.00){$Vol(\pi)$}
%%%%%%%%partition 1%%%%%%%%%%%%%%%%
\put(15.00,130.00){\line(0,1){12.00}}
\put(19.00,130.00){\line(0,1){8.00}}
\put(23.00,130.00){\line(0,1){4.00}}
\put(27.00,130.00){\line(0,1){4.00}}
\put(31.00,130.00){\line(0,1){8.00}}
\put(35.00,130.00){\line(0,1){12.00}}
\put(15.00,142.00){\line(1,0){20.00}}
\put(19.00,138.00){\line(1,0){12.00}}
\put(23.00,134.00){\line(1,0){4.00}}

\put(14.00,128.00){*}
\put(22.00,128.00){*}
\put(30.00,128.00){*}

\put(11.00,130.00){\line(0,1){1.5}}
\put(11.00,133.00){\line(0,1){1.5}}
\put(11.00,136.00){\line(0,1){1.5}}
\put(11.00,139.00){\line(0,1){1.5}}
\put(11.00,142.00){\line(0,1){1.5}}
\put(11.00,145.00){\line(0,1){1.0}}
\put(39.00,130.00){\line(0,1){1.5}}
\put(39.00,133.00){\line(0,1){1.5}}
\put(39.00,136.00){\line(0,1){1.5}}
\put(39.00,139.00){\line(0,1){1.5}}
\put(39.00,142.00){\line(0,1){1.5}}
\put(39.00,145.00){\line(0,1){1.0}}

\put(11.00,146.00){\line(1,0){1.5}}
\put(14.00,146.00){\line(1,0){1.5}}
\put(17.00,146.00){\line(1,0){1.5}}
\put(20.00,146.00){\line(1,0){1.5}}
\put(23.00,146.00){\line(1,0){1.5}}
\put(26.00,146.00){\line(1,0){1.5}}
\put(29.00,146.00){\line(1,0){1.5}}
\put(32.00,146.00){\line(1,0){1.5}}
\put(35.00,146.00){\line(1,0){1.5}}
\put(38.00,146.00){\line(1,0){1.0}}
\put(1.00,136.00){${\pi_{1}:}$}

\put(33.00,147.50){$\scriptstyle{x_1}$}
\put(30.00,143.50){$\scriptstyle{x_2}$}
\put(27.00,139.50){$\scriptstyle{x_3}$}
\put(24.00,135.50){$\scriptstyle{x_4}$}

\put(53.00,138.00){${x_2<x_1}$}
\put(53.00,134.00){${x_2<x_3}$}
\put(53.00,130.00){${x_4<x_3}$}

\put(80.00,146.00){${x_4<x_2<x_3<x_1}$}
\put(80.00,142.00){${x_2<x_4<x_3<x_1}$}
\put(80.00,138.00){${x_2<x_4<x_1<x_3}$}
\put(80.00,134.00){${x_4<x_2<x_1<x_3}$}
\put(80.00,130.00){${x_2<x_1<x_4<x_3}$}

\put(130.00,132.00){${\frac{5}{24}}$} 

%%%%%%%%partition 2%%%%%%%%%%%%%%%%
\put(15.00,100.00){\line(0,1){8.00}}
\put(19.00,100.00){\line(0,1){4.00}}
\put(23.00,100.00){\line(0,1){4.00}}
\put(27.00,100.00){\line(0,1){4.00}}
\put(31.00,100.00){\line(0,1){4.00}}
\put(35.00,100.00){\line(0,1){8.00}}
\put(15.00,108.00){\line(1,0){20.00}}
\put(19.00,104.00){\line(1,0){4.00}}
\put(27.00,104.00){\line(1,0){4.00}}

\put(14.00,98.00){*}
\put(22.00,98.00){*}
\put(30.00,98.00){*}

\put(11.00,100.00){\line(0,1){1.5}}
\put(11.00,103.00){\line(0,1){1.5}}
\put(11.00,106.00){\line(0,1){1.5}}
\put(11.00,109.00){\line(0,1){1.0}}
\put(11.00,111.00){\line(0,1){1.0}}
\put(39.00,100.00){\line(0,1){1.5}}
\put(39.00,103.00){\line(0,1){1.5}}
\put(39.00,106.00){\line(0,1){1.5}}
\put(39.00,109.00){\line(0,1){1.0}}
\put(39.00,111.00){\line(0,1){1.0}}

\put(11.00,112.00){\line(1,0){1.5}}
\put(14.00,112.00){\line(1,0){1.5}}
\put(17.00,112.00){\line(1,0){1.5}}
\put(20.00,112.00){\line(1,0){1.5}}
\put(23.00,112.00){\line(1,0){1.5}}
\put(26.00,112.00){\line(1,0){1.5}}
\put(29.00,112.00){\line(1,0){1.5}}
\put(32.00,112.00){\line(1,0){1.5}}
\put(35.00,112.00){\line(1,0){1.5}}
\put(38.00,112.00){\line(1,0){1.0}}
\put(1.00,105.00){${\pi_{2}:}$}

\put(32.00,113.30){$\scriptstyle{x_1}$}
\put(30.00,109.30){$\scriptstyle{x_2}$}
\put(28.00,105.30){$\scriptstyle{x_3}$}
\put(20.00,105.30){$\scriptstyle{x_4}$}

\put(53.00,108.00){${x_2<x_1}$}
\put(53.00,104.00){${x_2<x_4}$}
\put(53.00,100.00){${x_2<x_3}$}

\put(80.00,120.00){${x_2<x_3<x_4<x_1}$}
\put(80.00,116.00){${x_2<x_4<x_3<x_1}$}
\put(80.00,112.00){${x_2<x_4<x_1<x_3}$}
\put(80.00,108.00){${x_2<x_3<x_1<x_4}$}
\put(80.00,104.00){${x_2<x_1<x_3<x_4}$}
\put(80.00,100.00){${x_2<x_1<x_4<x_3}$}

\put(130.00,102.00){${\frac{6}{24}}$} 
%%%%%%%%partition 3%%%%%%%%%%%%%%%%
\put(15.00,74.00){\line(0,1){8.00}}
\put(19.00,74.00){\line(0,1){4.00}}
\put(23.00,74.00){\line(0,1){4.00}}
\put(27.00,74.00){\line(0,1){8.00}}
\put(31.00,74.00){\line(0,1){8.00}}
\put(35.00,74.00){\line(0,1){8.00}}
\put(15.00,82.00){\line(1,0){12.00}}
\put(19.00,78.00){\line(1,0){4.00}}
\put(31.00,82.00){\line(1,0){4.00}}

\put(14.00,72.00){*}
\put(22.00,72.00){*}
\put(30.00,72.00){*}

\put(11.00,74.00){\line(0,1){1.5}}
\put(11.00,77.00){\line(0,1){1.5}}
\put(11.00,80.00){\line(0,1){1.5}}
\put(11.00,83.00){\line(0,1){1.5}}
\put(11.00,85.00){\line(0,1){1.0}}
\put(39.00,74.00){\line(0,1){1.5}}
\put(39.00,77.00){\line(0,1){1.5}}
\put(39.00,80.00){\line(0,1){1.5}}
\put(39.00,83.00){\line(0,1){1.5}}
\put(39.00,85.00){\line(0,1){1.0}}

\put(11.00,86.00){\line(1,0){1.5}}
\put(14.00,86.00){\line(1,0){1.5}}
\put(17.00,86.00){\line(1,0){1.5}}
\put(20.00,86.00){\line(1,0){1.5}}
\put(23.00,86.00){\line(1,0){1.5}}
\put(26.00,86.00){\line(1,0){1.5}}
\put(29.00,86.00){\line(1,0){1.5}}
\put(32.00,86.00){\line(1,0){1.5}}
\put(35.00,86.00){\line(1,0){1.5}}
\put(38.00,86.00){\line(1,0){1.0}}
\put(1.00,79.00){${\pi_{3}:}$}

\put(33.00,87.30){$\scriptstyle{x_1}$}
\put(32.00,83.30){$\scriptstyle{x_2}$}
\put(24.00,83.30){$\scriptstyle{x_3}$}
\put(20.00,79.30){$\scriptstyle{x_4}$}

\put(53.00,82.00){${x_2<x_1}$}
\put(53.00,78.00){${x_3<x_4}$}
\put(53.00,74.00){${x_3<x_1}$}

\put(80.00,90.00){${x_2<x_3<x_4<x_1}$}
\put(80.00,86.00){${x_3<x_2<x_4<x_1}$}
\put(80.00,82.00){${x_3<x_4<x_2<x_1}$}
\put(80.00,78.00){${x_3<x_2<x_1<x_4}$}
\put(80.00,74.00){${x_2<x_3<x_1<x_4}$}

\put(130.00,76.00){${\frac{5}{24}}$} 
%%%%%%%%partition 4%%%%%%%%%%%%%%%%
\put(15.00,48.00){\line(0,1){8.00}}
\put(19.00,48.00){\line(0,1){8.00}}
\put(23.00,48.00){\line(0,1){8.00}}
\put(27.00,48.00){\line(0,1){4.00}}
\put(31.00,48.00){\line(0,1){4.00}}
\put(35.00,48.00){\line(0,1){8.00}}
\put(15.00,56.00){\line(1,0){4.00}}
\put(23.00,56.00){\line(1,0){12.00}}
\put(27.00,52.00){\line(1,0){4.00}}

\put(14.00,46.00){*}
\put(22.00,46.00){*}
\put(30.00,46.00){*}

\put(11.00,48.00){\line(0,1){1.5}}
\put(11.00,51.00){\line(0,1){1.5}}
\put(11.00,54.00){\line(0,1){1.5}}
\put(11.00,57.00){\line(0,1){1.5}}
\put(11.00,59.00){\line(0,1){1.0}}

\put(39.00,48.00){\line(0,1){1.5}}
\put(39.00,51.00){\line(0,1){1.5}}
\put(39.00,54.00){\line(0,1){1.5}}
\put(39.00,57.00){\line(0,1){1.5}}
\put(39.00,59.00){\line(0,1){1.0}}

\put(11.00,60.00){\line(1,0){1.5}}
\put(14.00,60.00){\line(1,0){1.5}}
\put(17.00,60.00){\line(1,0){1.5}}
\put(20.00,60.00){\line(1,0){1.5}}
\put(23.00,60.00){\line(1,0){1.5}}
\put(26.00,60.00){\line(1,0){1.5}}
\put(29.00,60.00){\line(1,0){1.5}}
\put(32.00,60.00){\line(1,0){1.5}}
\put(35.00,60.00){\line(1,0){1.5}}
\put(38.00,60.00){\line(1,0){1.0}}
\put(1.00,52.00){${\pi_{4}:}$}

\put(34.00,61.30){$\scriptstyle{x_1}$}
\put(31.00,57.30){$\scriptstyle{x_2}$}
\put(28.00,53.30){$\scriptstyle{x_3}$}
\put(16.00,57.30){$\scriptstyle{x_4}$}

\put(53.00,56.00){${x_2<x_1}$}
\put(53.00,52.00){${x_4<x_1}$}
\put(53.00,48.00){${x_2<x_3}$}

\put(80.00,64.00){${x_4<x_2<x_3<x_1}$}
\put(80.00,60.00){${x_2<x_4<x_3<x_1}$}
\put(80.00,56.00){${x_2<x_3<x_4<x_1}$}
\put(80.00,52.00){${x_4<x_2<x_1<x_3}$}
\put(80.00,48.00){${x_2<x_4<x_1<x_3}$}

\put(130.00,50.00){${\frac{5}{24}}$} 
%%%%%%%%partition 5%%%%%%%%%%%%%%%%
\put(15.00,18.00){\line(0,1){4.00}}
\put(19.00,18.00){\line(0,1){4.00}}
\put(23.00,18.00){\line(0,1){4.00}}
\put(27.00,18.00){\line(0,1){4.00}}
\put(31.00,18.00){\line(0,1){4.00}}
\put(35.00,18.00){\line(0,1){4.00}}
\put(15.00,22.00){\line(1,0){4.00}}
\put(23.00,22.00){\line(1,0){4.00}}
\put(31.00,22.00){\line(1,0){4.00}}

\put(14.00,16.00){*}
\put(22.00,16.00){*}
\put(30.00,16.00){*}

\put(11.00,18.00){\line(0,1){1.0}}
\put(11.00,20.00){\line(0,1){1.0}}
\put(11.00,22.00){\line(0,1){1.0}}
\put(11.00,24.50){\line(0,1){1.5}}
\put(39.00,18.00){\line(0,1){1.0}}
\put(39.00,20.00){\line(0,1){1.0}}
\put(39.00,22.00){\line(0,1){1.0}}
\put(39.00,24.50){\line(0,1){1.5}}

\put(11.00,26.00){\line(1,0){1.5}}
\put(14.00,26.00){\line(1,0){1.5}}
\put(17.00,26.00){\line(1,0){1.5}}
\put(20.00,26.00){\line(1,0){1.5}}
\put(23.00,26.00){\line(1,0){1.5}}
\put(26.00,26.00){\line(1,0){1.5}}
\put(29.00,26.00){\line(1,0){1.5}}
\put(32.00,26.00){\line(1,0){1.5}}
\put(35.00,26.00){\line(1,0){1.5}}
\put(38.00,26.00){\line(1,0){1.0}}
\put(1.00,21.00){${\pi_{5}:}$}

\put(34.00,27.30){$\scriptstyle{x_1}$}
\put(32.00,23.30){$\scriptstyle{x_2}$}
\put(24.00,23.30){$\scriptstyle{x_3}$}
\put(16.00,23.30){$\scriptstyle{x_4}$}

\put(53.00,26.00){${x_2<x_1}$}
\put(53.00,22.00){${x_3<x_1}$}
\put(53.00,18.00){${x_4<x_1}$}

\put(80.00,38.00){${x_4<x_3<x_2<x_1}$}
\put(80.00,34.00){${x_3<x_4<x_2<x_1}$}
\put(80.00,30.00){${x_4<x_2<x_3<x_1}$}
\put(80.00,26.00){${x_2<x_4<x_3<x_1}$}
\put(80.00,22.00){${x_3<x_2<x_4<x_1}$}
\put(80.00,18.00){${x_2<x_3<x_4<x_1}$}

\put(130.00,20.00){${\frac{6}{24}}$} 

\end{picture}
\caption{Noncrossing pair partitions of $[6]$, each extended by an imaginary block, corresponding
to $\varphi((T^*T)^3)$. Starred legs correspond to $T^*$. We assign continuous colors $x_j$ to blocks, where $j\in [4]$. The 
corresponding regions $R(\pi)$ inside the cube $[0,1]^{4}$ are given by three inequalities for colors of the blocks of $\pi$
and have volumes $Vol(\pi)$. Each region consists of simplices $S(\pi)$ defined by linearly ordered colors. Altogether we get $27=3^3$ simplices.}
\end{figure}

\begin{Example}
{\rm 
The set $\mathcal{O}_{4}$ of ordered rooted trees on 4 vertices consists of the 
trees given in Fig. 3, where the root is distinguished by a larger circle. 
In ordered trees, the children of any vertex are ordered and that is why $T_3$ and $T_4$ are 
inequivalent since different children of the roots have off-springs.
In turn, $\mathcal{NC}_{6}^{2}$ consists of the non-crossing pair partitions shown in Fig. 4,
where we also draw the imaginary blocks which can be identified with the roots of the corresponding trees. 
The natural bijection $\gamma:\mathcal{O}_{4}\rightarrow {\mathcal{NC}}_{6}^{2}$ is given by $\gamma(T_{k})=\pi_{k}$. In fact, that is why there are 5 trees of this type since $C_3=5$ and Catalan numbers count ordered rooted trees due to the 
bijection $\gamma$. In turn, there are $6!/3!=120$ different labeled ordered trees on 4 vertices if we label them by the 
4-element set in an arbitrary way.}
\end{Example}

\begin{figure}
\unitlength=1mm
\special{em.linewidth 0.5pt}
\linethickness{0.5pt}
\begin{picture}(140.00,110.00)(-15.00,-55.00)
%tree 1%%%%%%%%%%%%%%%%%%%%%%%%%%%%%
\put(5.00,30.00){\line(0,1){15.00}}
\put(5.00,30.00){\circle*{1.00}}
\put(5.00,35.00){\circle*{1.00}}
\put(5.00,40.00){\circle*{1.00}}
\put(5.00,45.00){\circle*{1.75}}
\put(2.00,30.00){$\scriptstyle{4}$}
\put(2.00,35.00){$\scriptstyle{2}$}
\put(2.00,40.00){$\scriptstyle{3}$}
\put(4.50,47.00){$\scriptstyle{1}$}
%tree 2%%%%%%%%%%%%%%%%%%%%%%%%%%%%%
\put(25.00,30.00){\line(0,1){15.00}}
\put(25.00,30.00){\circle*{1.00}}
\put(25.00,35.00){\circle*{1.00}}
\put(25.00,40.00){\circle*{1.00}}
\put(25.00,45.00){\circle*{1.75}}
\put(22.00,30.00){$\scriptstyle{3}$}
\put(22.00,35.00){$\scriptstyle{2}$}
\put(22.00,40.00){$\scriptstyle{4}$}
\put(24.50,47.00){$\scriptstyle{1}$}
%tree 3%%%%%%%%%%%%%%%%%%%%%%%%%%%%%
\put(45.00,30.00){\line(0,1){15.00}}
\put(45.00,30.00){\circle*{1.00}}
\put(45.00,35.00){\circle*{1.00}}
\put(45.00,40.00){\circle*{1.00}}
\put(45.00,45.00){\circle*{1.75}}
\put(42.00,30.00){$\scriptstyle{3}$}
\put(42.00,35.00){$\scriptstyle{1}$}
\put(42.00,40.00){$\scriptstyle{4}$}
\put(44.50,47.00){$\scriptstyle{2}$}
%tree 4%%%%%%%%%%%%%%%%%%%%%%%%%%%%%
\put(65.00,30.00){\line(0,1){15.00}}
\put(65.00,30.00){\circle*{1.00}}
\put(65.00,35.00){\circle*{1.00}}
\put(65.00,40.00){\circle*{1.00}}
\put(65.00,45.00){\circle*{1.75}}
\put(62.00,30.00){$\scriptstyle{4}$}
\put(62.00,35.00){$\scriptstyle{1}$}
\put(62.00,40.00){$\scriptstyle{3}$}
\put(64.50,47.00){$\scriptstyle{2}$}
%tree 5%%%%%%%%%%%%%%%%%%%%%%%%%%%%%
\put(85.00,30.00){\line(0,1){15.00}}
\put(85.00,30.00){\circle*{1.00}}
\put(85.00,35.00){\circle*{1.00}}
\put(85.00,40.00){\circle*{1.00}}
\put(85.00,45.00){\circle*{1.75}}
\put(82.00,30.00){$\scriptstyle{2}$}
\put(82.00,35.00){$\scriptstyle{1}$}
\put(82.00,40.00){$\scriptstyle{4}$}
\put(84.50,47.00){$\scriptstyle{3}$}
%tree 1%%%%%%%%%%%%%%%%%%%%%%%%%%%%%
\put(5.00,15.00){\line(0,1){5.00}}
\put(5.00,15.00){\line(1,-1){5.00}}
\put(5.00,15.00){\line(-1,-1){5.00}}
\put(10.00,10.00){\circle*{1.00}}
\put(0.00,10.00){\circle*{1.00}}
\put(5.00,15.00){\circle*{1.00}}
\put(5.00,20.00){\circle*{1.75}}
\put(9.50,7.00){$\scriptstyle{3}$}
\put(-0.50,7.00){$\scriptstyle{2}$}
\put(2.50,15.00){$\scriptstyle{4}$}
\put(4.50,22.00){$\scriptstyle{1}$}
%tree 2%%%%%%%%%%%%%%%%%%%%%%%%%%%%%
\put(25.00,15.00){\line(0,1){5.00}}
\put(25.00,15.00){\line(1,-1){5.00}}
\put(25.00,15.00){\line(-1,-1){5.00}}
\put(30.00,10.00){\circle*{1.00}}
\put(20.00,10.00){\circle*{1.00}}
\put(25.00,15.00){\circle*{1.00}}
\put(25.00,20.00){\circle*{1.75}}
\put(29.50,7.00){$\scriptstyle{2}$}
\put(19.50,7.00){$\scriptstyle{3}$}
\put(22.50,15.00){$\scriptstyle{4}$}
\put(24.50,22.00){$\scriptstyle{1}$}
%tree 3%%%%%%%%%%%%%%%%%%%%%%%%%%%%%
\put(45.00,15.00){\line(0,1){5.00}}
\put(45.00,15.00){\line(1,-1){5.00}}
\put(45.00,15.00){\line(-1,-1){5.00}}
\put(50.00,10.00){\circle*{1.00}}
\put(40.00,10.00){\circle*{1.00}}
\put(45.00,15.00){\circle*{1.00}}
\put(45.00,20.00){\circle*{1.75}}
\put(49.50,7.00){$\scriptstyle{1}$}
\put(39.50,7.00){$\scriptstyle{3}$}
\put(42.50,15.00){$\scriptstyle{4}$}
\put(44.50,22.00){$\scriptstyle{2}$}
%tree 4%%%%%%%%%%%%%%%%%%%%%%%%%%%%%
\put(65.00,15.00){\line(0,1){5.00}}
\put(65.00,15.00){\line(1,-1){5.00}}
\put(65.00,15.00){\line(-1,-1){5.00}}
\put(70.00,10.00){\circle*{1.00}}
\put(60.00,10.00){\circle*{1.00}}
\put(65.00,15.00){\circle*{1.00}}
\put(65.00,20.00){\circle*{1.75}}
\put(69.50,7.00){$\scriptstyle{3}$}
\put(59.50,7.00){$\scriptstyle{1}$}
\put(62.50,15.00){$\scriptstyle{4}$}
\put(64.50,22.00){$\scriptstyle{2}$}
%tree 5%%%%%%%%%%%%%%%%%%%%%%%%%%%%%
\put(85.00,15.00){\line(0,1){5.00}}
\put(85.00,15.00){\line(1,-1){5.00}}
\put(85.00,15.00){\line(-1,-1){5.00}}
\put(90.00,10.00){\circle*{1.00}}
\put(80.00,10.00){\circle*{1.00}}
\put(85.00,15.00){\circle*{1.00}}
\put(85.00,20.00){\circle*{1.75}}
\put(89.50,7.00){$\scriptstyle{2}$}
\put(79.50,7.00){$\scriptstyle{1}$}
\put(82.50,15.00){$\scriptstyle{4}$}
\put(84.50,22.00){$\scriptstyle{3}$}
%tree 6%%%%%%%%%%%%%%%%%%%%%%%%%%%%%
\put(105.00,15.00){\line(0,1){5.00}}
\put(105.00,15.00){\line(1,-1){5.00}}
\put(105.00,15.00){\line(-1,-1){5.00}}
\put(110.00,10.00){\circle*{1.00}}
\put(100.00,10.00){\circle*{1.00}}
\put(105.00,15.00){\circle*{1.00}}
\put(105.00,20.00){\circle*{1.75}}
\put(109.50,7.00){$\scriptstyle{1}$}
\put(99.50,7.00){$\scriptstyle{2}$}
\put(102.50,15.00){$\scriptstyle{4}$}
\put(104.50,22.00){$\scriptstyle{3}$}
%tree 1%%%%%%%%%%%%%%%%%%%%%%%%%%%%%%%%
\put(5.00,0.00){\line(1,-1){5.00}}
\put(5.00,0.00){\line(-1,-1){5.00}}
\put(0.00,-5.00){\line(0,-1){5.00}}
\put(0.00,-5.00){\circle*{1.00}}
\put(0.00,-10.00){\circle*{1.00}}
\put(10.00,-5.00){\circle*{1.00}}
\put(5.00,0.00){\circle*{1.75}}
\put(-2.00,-5.00){$\scriptstyle{3}$}
\put(-0.50,-13.00){$\scriptstyle{2}$}
\put(11.00,-5.00){$\scriptstyle{4}$}
\put(4.50,2.00){$\scriptstyle{1}$}
%tree 2%%%%%%%%%%%%%%%%%%%%%%%%%%%%%%%%
\put(25.00,0.00){\line(1,-1){5.00}}
\put(25.00,0.00){\line(-1,-1){5.00}}
\put(20.00,-5.00){\line(0,-1){5.00}}
\put(20.00,-5.00){\circle*{1.00}}
\put(20.00,-10.00){\circle*{1.00}}
\put(30.00,-5.00){\circle*{1.00}}
\put(25.00,0.00){\circle*{1.75}}
\put(18.00,-5.00){$\scriptstyle{4}$}
\put(19.50,-13.00){$\scriptstyle{2}$}
\put(31.00,-5.00){$\scriptstyle{3}$}
\put(24.50,2.00){$\scriptstyle{1}$}
%tree 3%%%%%%%%%%%%%%%%%%%%%%%%%%%%%%%%
\put(45.00,0.00){\line(1,-1){5.00}}
\put(45.00,0.00){\line(-1,-1){5.00}}
\put(40.00,-5.00){\line(0,-1){5.00}}
\put(40.00,-5.00){\circle*{1.00}}
\put(40.00,-10.00){\circle*{1.00}}
\put(50.00,-5.00){\circle*{1.00}}
\put(45.00,0.00){\circle*{1.75}}
\put(38.00,-5.00){$\scriptstyle{4}$}
\put(39.50,-13.00){$\scriptstyle{3}$}
\put(51.00,-5.00){$\scriptstyle{2}$}
\put(44.50,2.00){$\scriptstyle{1}$}
%tree 4%%%%%%%%%%%%%%%%%%%%%%%%%%%%%%%%
\put(65.00,0.00){\line(1,-1){5.00}}
\put(65.00,0.00){\line(-1,-1){5.00}}
\put(60.00,-5.00){\line(0,-1){5.00}}
\put(60.00,-5.00){\circle*{1.00}}
\put(60.00,-10.00){\circle*{1.00}}
\put(70.00,-5.00){\circle*{1.00}}
\put(65.00,0.00){\circle*{1.75}}
\put(58.00,-5.00){$\scriptstyle{4}$}
\put(59.50,-13.00){$\scriptstyle{1}$}
\put(71.00,-5.00){$\scriptstyle{3}$}
\put(64.50,2.00){$\scriptstyle{2}$}
%tree 5%%%%%%%%%%%%%%%%%%%%%%%%%%%%%%%%
\put(85.00,0.00){\line(1,-1){5.00}}
\put(85.00,0.00){\line(-1,-1){5.00}}
\put(80.00,-5.00){\line(0,-1){5.00}}
\put(80.00,-5.00){\circle*{1.00}}
\put(80.00,-10.00){\circle*{1.00}}
\put(90.00,-5.00){\circle*{1.00}}
\put(85.00,0.00){\circle*{1.75}}
\put(78.00,-5.00){$\scriptstyle{3}$}
\put(79.50,-13.00){$\scriptstyle{1}$}
\put(91.00,-5.00){$\scriptstyle{4}$}
\put(84.50,2.00){$\scriptstyle{2}$}
%tree 1%%%%%%%%%%%%%%%%%%%%%%%%%%%%%%%%
\put(5.00,-20.00){\line(1,-1){5.00}}
\put(5.00,-20.00){\line(-1,-1){5.00}}
\put(10.00,-25.00){\line(0,-1){5.00}}
\put(10.00,-25.00){\circle*{1.00}}
\put(0.00,-25.00){\circle*{1.00}}
\put(10.00,-25.00){\circle*{1.00}}
\put(5.00,-20.00){\circle*{1.75}}
\put(-2.00,-25.00){$\scriptstyle{4}$}
\put(9.50,-33.00){$\scriptstyle{2}$}
\put(11.00,-25.00){$\scriptstyle{3}$}
\put(4.50,-18.00){$\scriptstyle{1}$}
%tree 2%%%%%%%%%%%%%%%%%%%%%%%%%%%%%%%%
\put(25.00,-20.00){\line(1,-1){5.00}}
\put(25.00,-20.00){\line(-1,-1){5.00}}
\put(30.00,-25.00){\line(0,-1){5.00}}
\put(20.00,-25.00){\circle*{1.00}}
\put(30.00,-30.00){\circle*{1.00}}
\put(30.00,-25.00){\circle*{1.00}}
\put(25.00,-20.00){\circle*{1.75}}
\put(18.00,-25.00){$\scriptstyle{3}$}
\put(29.50,-33.00){$\scriptstyle{2}$}
\put(31.00,-25.00){$\scriptstyle{4}$}
\put(24.50,-18.00){$\scriptstyle{1}$}
%tree 3%%%%%%%%%%%%%%%%%%%%%%%%%%%%%%%%
\put(45.00,-20.00){\line(1,-1){5.00}}
\put(45.00,-20.00){\line(-1,-1){5.00}}
\put(50.00,-25.00){\line(0,-1){5.00}}
\put(40.00,-25.00){\circle*{1.00}}
\put(50.00,-30.00){\circle*{1.00}}
\put(50.00,-25.00){\circle*{1.00}}
\put(45.00,-20.00){\circle*{1.75}}
\put(38.00,-25.00){$\scriptstyle{2}$}
\put(49.50,-33.00){$\scriptstyle{3}$}
\put(51.00,-25.00){$\scriptstyle{4}$}
\put(44.50,-18.00){$\scriptstyle{1}$}
%tree 4%%%%%%%%%%%%%%%%%%%%%%%%%%%%%%%%
\put(65.00,-20.00){\line(1,-1){5.00}}
\put(65.00,-20.00){\line(-1,-1){5.00}}
\put(70.00,-25.00){\line(0,-1){5.00}}
\put(60.00,-25.00){\circle*{1.00}}
\put(70.00,-30.00){\circle*{1.00}}
\put(70.00,-25.00){\circle*{1.00}}
\put(65.00,-20.00){\circle*{1.75}}
\put(58.00,-25.00){$\scriptstyle{4}$}
\put(69.50,-33.00){$\scriptstyle{1}$}
\put(71.00,-25.00){$\scriptstyle{3}$}
\put(64.50,-18.00){$\scriptstyle{2}$}
%tree 5%%%%%%%%%%%%%%%%%%%%%%%%%%%%%%%
\put(85.00,-20.00){\line(1,-1){5.00}}
\put(85.00,-20.00){\line(-1,-1){5.00}}
\put(90.00,-25.00){\line(0,-1){5.00}}
\put(80.00,-25.00){\circle*{1.00}}
\put(90.00,-30.00){\circle*{1.00}}
\put(90.00,-25.00){\circle*{1.00}}
\put(85.00,-20.00){\circle*{1.75}}
\put(78.00,-25.00){$\scriptstyle{3}$}
\put(89.50,-33.00){$\scriptstyle{1}$}
\put(91.00,-25.00){$\scriptstyle{4}$}
\put(84.50,-18.00){$\scriptstyle{2}$}

%tree 1%%%%%%%%%%%%%%%%%%%%%%%%%%%%%%%
\put(5.00,-40.00){\line(0,-1){5.00}}
\put(5.00,-40.00){\line(1,-1){5.00}}
\put(5.00,-40.00){\line(-1,-1){5.00}}
\put(0.00,-45.00){\circle*{1.00}}
\put(5.00,-45.00){\circle*{1.00}}
\put(10.00,-45.00){\circle*{1.00}}
\put(5.00,-40.00){\circle*{1.75}}
\put(4.50,-38.00){$\scriptstyle{1}$}
\put(-0.50,-48.00){$\scriptstyle{4}$}
\put(4.50,-48.00){$\scriptstyle{3}$}
\put(9.50,-48.00){$\scriptstyle{2}$}
%tree 2%%%%%%%%%%%%%%%%%%%%%%%%%%%%%%%
\put(25.00,-40.00){\line(0,-1){5.00}}
\put(25.00,-40.00){\line(1,-1){5.00}}
\put(25.00,-40.00){\line(-1,-1){5.00}}
\put(20.00,-45.00){\circle*{1.00}}
\put(25.00,-45.00){\circle*{1.00}}
\put(30.00,-45.00){\circle*{1.00}}
\put(25.00,-40.00){\circle*{1.75}}
\put(24.50,-38.00){$\scriptstyle{1}$}
\put(19.50,-48.00){$\scriptstyle{3}$}
\put(24.50,-48.00){$\scriptstyle{4}$}
\put(29.50,-48.00){$\scriptstyle{2}$}
%tree 3%%%%%%%%%%%%%%%%%%%%%%%%%%%%%%%
\put(45.00,-40.00){\line(0,-1){5.00}}
\put(45.00,-40.00){\line(1,-1){5.00}}
\put(45.00,-40.00){\line(-1,-1){5.00}}
\put(40.00,-45.00){\circle*{1.00}}
\put(45.00,-45.00){\circle*{1.00}}
\put(50.00,-45.00){\circle*{1.00}}
\put(45.00,-40.00){\circle*{1.75}}
\put(44.50,-38.00){$\scriptstyle{1}$}
\put(39.50,-48.00){$\scriptstyle{4}$}
\put(44.50,-48.00){$\scriptstyle{2}$}
\put(49.50,-48.00){$\scriptstyle{3}$}
%tree 4%%%%%%%%%%%%%%%%%%%%%%%%%%%%%%%
\put(65.00,-40.00){\line(0,-1){5.00}}
\put(65.00,-40.00){\line(1,-1){5.00}}
\put(65.00,-40.00){\line(-1,-1){5.00}}
\put(60.00,-45.00){\circle*{1.00}}
\put(65.00,-45.00){\circle*{1.00}}
\put(70.00,-45.00){\circle*{1.00}}
\put(65.00,-40.00){\circle*{1.75}}
\put(64.50,-38.00){$\scriptstyle{1}$}
\put(59.50,-48.00){$\scriptstyle{3}$}
\put(64.50,-48.00){$\scriptstyle{2}$}
\put(69.50,-48.00){$\scriptstyle{4}$}
%tree 5%%%%%%%%%%%%%%%%%%%%%%%%%%%%%%%
\put(85.00,-40.00){\line(0,-1){5.00}}
\put(85.00,-40.00){\line(1,-1){5.00}}
\put(85.00,-40.00){\line(-1,-1){5.00}}
\put(80.00,-45.00){\circle*{1.00}}
\put(85.00,-45.00){\circle*{1.00}}
\put(90.00,-45.00){\circle*{1.00}}
\put(85.00,-40.00){\circle*{1.75}}
\put(84.50,-38.00){$\scriptstyle{1}$}
\put(79.50,-48.00){$\scriptstyle{2}$}
\put(84.50,-48.00){$\scriptstyle{4}$}
\put(89.50,-48.00){$\scriptstyle{3}$}
%tree 6%%%%%%%%%%%%%%%%%%%%%%%%%%%%%%%
\put(105.00,-40.00){\line(0,-1){5.00}}
\put(105.00,-40.00){\line(1,-1){5.00}}
\put(105.00,-40.00){\line(-1,-1){5.00}}
\put(100.00,-45.00){\circle*{1.00}}
\put(105.00,-45.00){\circle*{1.00}}
\put(110.00,-45.00){\circle*{1.00}}
\put(105.00,-40.00){\circle*{1.75}}
\put(104.50,-38.00){$\scriptstyle{1}$}
\put(99.50,-48.00){$\scriptstyle{2}$}
\put(104.50,-48.00){$\scriptstyle{3}$}
\put(109.50,-48.00){$\scriptstyle{4}$}
\end{picture}
\caption{The set $\mathpzc{A}_{4}$ of alternating ordered rooted trees of type I. This set is in bijection with the set of noncrossing colored pair partitions with totally ordered colors given in Fig.~4, or with the corresponding simplices. The alternating ordered rooted trees are listed in the same order as the corresponding simplices in Fig.~4.}
\end{figure}
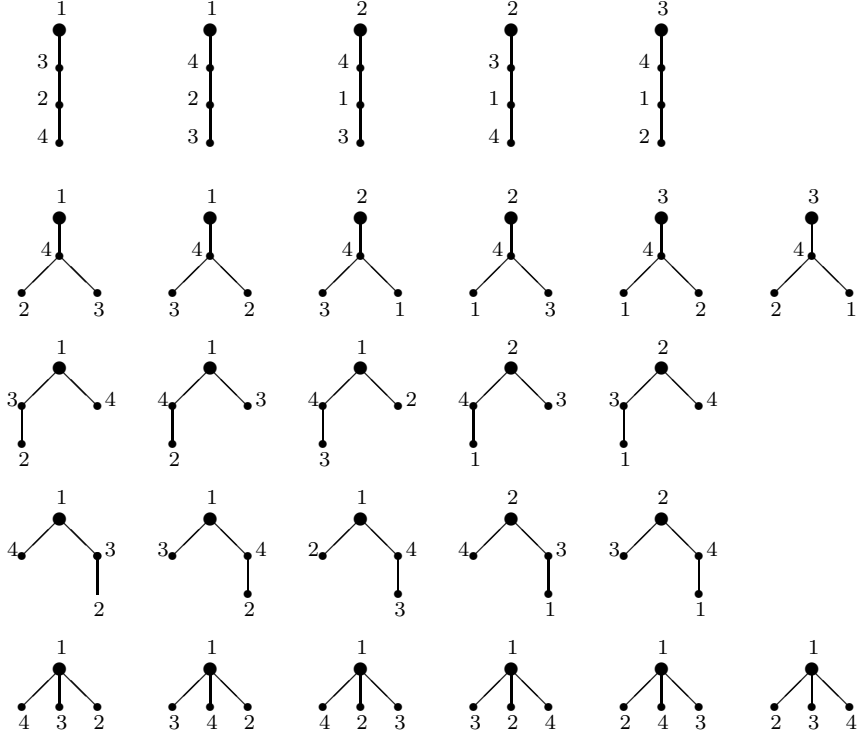

\begin{Definition}
{\rm 
Let $(v_1,v_2,v_3, v_4, \ldots)$ be a path in an ordered labeled rooted tree $T$ on 
$n+1$ vertices, which means that $v_j$ is a son of $v_{j-1}$, and let 
$(x_1,x_2,x_3,x_4, \ldots)$ be the corresponding sequence of labels.
Then $T$ is called {\it alternating} if this sequence satisfies one of the inequalities,
$$
x_1>x_2<x_3>x_4\ldots,\;\;\;\; {\rm or}\;\; \;\;x_1<x_2>x_3<x_4\ldots,
$$
i.e. the differences of labels corresponding to the neighboring vertices alternate in sign. 
These two types of conditions split the set of alternating ordered rooted trees on $n+1$ vertices into two subsets
of the same cardinalities, which we will call {\it of type I} and {\it of type II}, respectively. 
Denote by $\mathcal{A}_{n}$ the set of alternating ordered rooted trees on $n+1$ vertices of type I.}
\end{Definition}

\begin{Remark}
{\rm 
There is a nice enumeration result of Chauve, Dulucq and Rechnitzer [2] which says that 
$$
|\mathcal{A}_{n}|=n^{n}
$$
for any natural $n$. Let us recall that in our notation $\mathcal{A}_{n}$ stands for
the set of alternating ordered rooted trees of type I on $n+1$ vertices. Thus, 
we need to multiply this formula by 2 to get the number of all alternating ordered rooted trees. 
Note also that a typical formula refers to trees on $n$ vertices.
}
\end{Remark}

\begin{Example}
{\rm 
Among the ordered labeled rooted trees on 4 vertices, there are $2\times 3^{3}=54$ alternating ones. 
In this example, each type contains $3^{3}=27$ alternating ordered rooted trees. The complete set $\mathcal{A}_{4}$ of all alternating ordered rooted trees of type I is given in Figure 5. In view of the above bijection results, the cardinality of all non-crossing pair partitions 
of the 6-element set with alternating colorings is also 54 and there are $27$ partitions in which each block of odd depth has a smaller color than its nearest outer block (equivalently, the color of the imaginary block, which is assumed to have zero depth, is greater
than the colors of all blocks for which the imaginary block is the nearest outer block).}
\end{Example}

\begin{Corollary}
The moments of $T^*T$, where $T$ is the triangular operator, are 
$$
M_{n}=\varphi((T^{*}T)^{n})= \frac{n^{n}}{(n+1)!}
$$
for any $n\in {\mathbb N}$.
\end{Corollary}
{\it Proof.}
Let $\epsilon_{j}=*$ if $j$ is odd and $\epsilon_j=1$ if $j$ is even. In this special case, it is easy to see that
$\mathcal{NC}^{2}(\epsilon_1, \ldots, \epsilon_{2n})\cong \mathcal{NC}^{2}_{2n}$.
Observe that in the case of alternating starred and unstarred legs
there can be no pairs $(V,o(V))$ of type 3 and 4 since otherwise there would be unequal numbers of starred and unstarred legs between 
the right leg of $V$ and the left leg of $o(V)$, which would mean that there must be a block between $V$ and $o(V)$, 
which is a contradiction. Therefore, blocks with starred left and right legs must alternate as we take a sequence of blocks
$(V_{i_1}, \ldots , V_{i_p})$, where each block is the nearest outer block of its successor. By Corollary 6.1, 
we need to compute $Vol(\pi)$ for each $\pi\in \mathcal{NC}^{2}_{2n}$. Each of the corresponding regions $R(\pi)$ 
is defined by a set of $n$ inequalities for colors $x_1, \ldots , x_{n+1}$. Irrespective of what symbols represent 
these colors, in order to satisfy these inequalities, we have to find the number of total orderings 
of the form 
$$
x_{j_1}<x_{j_2}<\ldots <x_{j_{n+1}}
$$
which satisfy the given inequalities. The number of these total orderings is 
equal to the number of $n+1$-simplices, each of volume $1/(n+1)!$, into which $R(\pi)$ decomposes. 
The key observation is that in order to compute the number of these total orderings
corresponding to $\pi$ (under conditions given by $n$ inequalities which express orders between the colors of 
each $V$ and $o(V)$) it suffices to count in how many ways we can label blocks of $\pi$ with natural numbers from $[n+1]$ 
in such a way that orders between these numbers alternate as we go down each sequence of type $(V_{i_1}, \ldots , V_{i_p})$. 
More specifically, to colors $x_{j_1}, x_{j_2},\ldots, x_{j_{n+1}}$ in the total ordering (defining a simplex) given above we assign 
numbers $n+1, n, \ldots , 1$, respectively, thus we assign number 1 to the biggest color and the number $n+1$ - to the smallest one.
Now, if we use the bijection between ordered labeled rooted trees and colored noncrossing pair partitions, 
it suffices to enumerate all alternating ordered rooted trees on $n+1$ vertices. 
The enumeration result of [2] mentioned above completes the proof.
\hfill $\blacksquare$


\begin{thebibliography}{99}
\bibitem{[1]}
F. Benaych-Georges, Rectangular random matrices, related convolution, {\it Probab. Theory Relat. Fields} {\bf 144} (2009), 471-515.
\bibitem{[2]}
C. Chauve, S. Dulucq, A. Rechnitzer, 
Enumerating alternating trees, {\it J. Combin. Theory Ser. A} {\bf 94} (2001), no. 1, 142-151.
\bibitem{[3]}
K. Dykema, On certain free product factors via an extended matrix model, {\it J. Funct. Anal.} {\bf 112} (1993), 31-60. 
\bibitem{[4]}
K. Dykema, U. Haagerup, DT-operators and decomposability of Voiculescu's circular operator, {\it Amer. J. Math.} {\bf 126} (2004), 121-189.
\bibitem{[5]}
K. Dykema, U. Haagerup, Invariant subspaces of the quasinilpotent DT-operator, {\it J. Funct. Anal.} {\bf 209}
(2004), 332-366.
\bibitem{[6]}
R.~V. Kadison, J.~R. Ringrose, {\it Fundamentals of the Theory of Operator Algebras, Vol. II. Advanced Theory.}, Graduate Studies in Mathematics, 16, American Mathematical Society, Providence, RI, 1997.
\bibitem{[7]}
R. Lenczewski, Matricially free random variables, {\it J. Funct. Anal.} {\bf 258} (2010), 4075-4121.
\bibitem{[8]}
R. Lenczewski, Asymptotic properties of random matrices and pseudomatrices, 
{\it Adv. Math.} {\bf 228} (2011), 2403-2440.
\bibitem{[9]}
R. Lenczewski, Limit distributions of random matrices, {\it Adv. Math.} {\bf 263} (2014), 253-320.
\bibitem{[10]}
R. Lenczewski, Matricial circular systems and random matrices, {\it Random Matrices Theory Appl.}
{\bf 5}, no.~4 (2016), 1650012.
\bibitem{[11]}
J. Mingo, R. Speicher, \emph{Free Probability and Random Matrices}, 
Fields Institute Monographs, Vol.~35, Springer-Verlag, New York, 2017.
\bibitem{[12]}
A. Nica, R. Speicher, {\it Lectures on the Combinatorics of Free Probability}, Cambridge University Press, Cambridge 2006.
\bibitem{[13]}
D. Shlyakhtenko, Random Gaussian band matrices and freeness with amalgamation, {\it Int. Math. Res. Notices}
{\bf 20} (1996), 1013-1025.
\bibitem{[14]}
P. \'{S}niady, Multinomial identities arising from free probability theory, {\it J. Combin. Theory Ser. A} {\bf 101} (2003), no. 1, 1-19.
\bibitem{[15]}
P. \'{S}niady, Generalized Cauchy identities, trees and multidimesnional 
Brownian motions. I. Bijective proof of generalized Cauchy identities. 
{\it  Electron. J. Combin.} {\bf 13} (2006), no. 1, Research Paper 62, 27pp.
\bibitem{[16]}
D. Voiculescu, Limit laws for random matrices and free products, {\it Invent. Math.} {\bf 104} (1991), 201-220.
\bibitem{[17]}
D. Voiculescu, Circular and semicircular systems and free product factors, {\it Progress in Math.} {\bf 92}, Birkhauser, 1990.
\bibitem{[18]}
D. Voiculescu, K. Dykema, A. Nica, {\it Free random variables}, CRM Monograph
Series, No.1, A.M.S., Providence, 1992.
\end{thebibliography}
\end{document}